\documentclass{article}
\usepackage{latexsym,amsfonts,amsmath,amsthm,makeidx}
\usepackage{makeidx}
\usepackage[top=1.0in, bottom=1.0in, left=1.25in, right=1.25in]{geometry}
\makeindex
\newtheorem{definition}{Definition}[section]
\newtheorem{theorem}{Theorem}[section]
\newtheorem{lemma}{Lemma}[section]
\newtheorem{corollary}{Corollary}[section]

\newtheorem{remark}{Remark}[section]
\newcommand{\RN}{\mathbb R^N}
\newcommand{\rw}{\rightarrow}
\newcommand{\Om}{\Omega}

\newcommand{\iy}{\infty}

\newcommand{\s}{\section}

\newcommand{\la}{\lambda}

\newcommand{\R}{\mathbb R}
\newcommand{\al}{\alpha}

\newcommand{\bt}{\begin{theorem}}
\newcommand{\et}{\end{theorem}}
\newcommand{\bl}{\begin{lemma}}
\newcommand{\el}{\end{lemma}}
\newcommand{\bd}{\begin{definition}}
\newcommand{\ed}{\end{definition}}
\newcommand{\bc}{\begin{corollary}}
\newcommand{\ec}{\end{corollary}}
\newcommand{\bp}{\begin{proof}}
\newcommand{\ep}{\end{proof}}
\newcommand{\bx}{\begin{example}}
\newcommand{\ex}{\end{example}}
\newcommand{\bi}{\begin{exercise}}
\newcommand{\ei}{\end{exercise}}
\newcommand{\bo}{\begin{prop}}
\newcommand{\eo}{\end{prop}}
\newcommand{\br}{\begin{remark}}
\newcommand{\er}{\end{remark}}
\newcommand{\be}{\begin{equation}}
\newcommand{\ee}{\end{equation}}
\newcommand{\ba}{\begin{align}}
\newcommand{\ea}{\end{align}}
\newcommand{\bn}{\begin{enumerate}}
\newcommand{\en}{\end{enumerate}}
\newcommand{\bg}{\begin{align*}}
\newcommand{\bcs}{\begin{cases}}
\newcommand{\ecs}{\end{cases}}

\newcommand{\ga}{\gamma}

\newcommand{\bean}{\begin{eqnarray*}}
\newcommand{\eean}{\end{eqnarray*}}

\hfuzz=\maxdimen
\numberwithin{equation}{section}

\begin{document}

\title{\bf On the existence and regularity of vector solutions for quasilinear systems with linear coupling\thanks {E-mail:
aoy15@mails.tsinghua.edu.cn (Ao);\;\; 15901026320@163.com(Wang); \quad wzou@math.tsinghua.edu.cn (Zou)}}
\date{}
\author{{\bf Yong Ao$^{1}$,  Jiaqi Wang$^{2}$,   Wenming Zou$^{3}$}\\
\footnotesize {\it $^{1, 2, 3}$Department of Mathematical Sciences, Tsinghua
University,}\\
\footnotesize {\it Beijing 100084, China}\\ }

\maketitle \vskip0.16in
\begin{center}
\begin{minipage}{120mm}
\begin{center}{\bf Abstract}\end{center}

We study the following coupled
system of quasilinear equations:
\begin{equation}\nonumber
\begin{cases}-\Delta_p u +|u|^{p-2}u =
f(u)+\lambda v, \quad x\in \R^N,\\
-\Delta_p v +|v|^{p-2}v =
g(v)+\lambda u, \quad x\in \R^N.
\end{cases}
\end{equation}
Under some assumptions on the nonlinear
terms $f$ and $g$, we establish some results about the existence and regularity of vector
solutions for the p-Laplacian systems by using variational methods. In particular, we get two pairs of nontrivial
solutions.  We also study their different  asymptotic behavior of solutions as the coupling parameter $\lambda$ tends to zero.\\

\vskip0.23in

{\it Key  words:}   p-Laplacian system, Least energy solutions, Moser iteration,  variational methods.

\vskip0.123in

\noindent {Mathematics Subject Classification 2010: 35B33, 35J20, 58E05}  
\end{minipage}
\end{center}

\vskip0.16in \s{Introduction}

In this paper we study the following coupled
system of quasilinear equations:
\be\label{1-1}
\begin{cases}-\Delta_p u +|u|^{p-2}u =
f(u)+\lambda v, \quad x\in \R^N,\\
-\Delta_p v +|v|^{p-2}v =
g(v)+\lambda u, \quad x\in \R^N.\end{cases}\ee
where the coupling constant $\la \ge 0$, $N \ge 3$, $1<p<N$ and $\Delta_pu=\mbox{div}(|\nabla u|^{p-2}\nabla u)$ is the p-Laplacian operator. For $p\neq 2$, the operator $\Delta_pu$ arises in non-Newtonian fluids, flow through porous media, nonlinear elasticity, and other physical phenomena. A solution $(u,v)\in W^{1,p}(\R^N)\times W^{1,p}(\R^N)\setminus \{(0,0)\}$ of system \eqref{1-1} is called a nontrivial solution, and a positive solution if $u>0,v>0$. A solution is called a ground state if $(u,v)\neq (0,0)$ and its energy is minimal among the energy of all the nontrivial solutions of \eqref{1-1}. Obviously, the solutions of \eqref{1-1} are the critical points of the functional $I_\la:W^{1,p}(\R^N)\times W^{1,p}(\R^N)\rightarrow \R$ given by
\be\label{1-2}
\begin{split}
I_\la(u,v):=&\frac{1}{p}\int_{\R^N}\big(|\nabla u|^p+|\nabla v|^p\big)dx+\frac{1}{p}\int_{\R^N}\big(|u|^p+|v|^p\big)dx\\
&-\int_{\R^N}\big(F(u)+G(v)+\la uv\big)dx,
\end{split}
\ee
if $I_\la\in C^1(W^{1,p}(\R^N)\times W^{1,p}(\R^N), \R)$. Here
$$F(u)=\int_{0}^{u}f(s)ds\;\;\mbox{and}\;\;G(v)=\int_{0}^{v}g(s)ds.$$

In recent years great interest has been devoted to the study of elliptic systems involving the p-Laplacian operator：
\be\label{1-3}
-\mbox{div}(|\nabla u_i|^{p-2}\nabla u_i)=g_i(u),~~~i=1,...,m,
\ee
where $u=(u_1,...,u_m):\R^N\longrightarrow \R^m,\;1<p\leq N$ and $g_i(u)=\frac{\partial G}{\partial u_i}(u)$ for some function $G\in C^1(\R^m)$. A series of papers have been devoted to the case $p = 2$ and fairly optimal conditions on $g_i$ have been found by Berestycki and Lions \cite{BeL} for $m = 1$ and by Br\'{e}zis and Lieb \cite{BrL} for $m \geq 1$. In \cite{BJM}, Byeon Jeanjean and Maris proved that the least energy solutions for \eqref{1-3} are radially symmetric under some assumptions on the corresponding minimizing problem. 

\vskip0.12in
When $p=2$, the system \eqref{1-1} turns to be the following Schr\"{o}dinger system:
\be\label{1-4}
\begin{cases}-\Delta u +u =
	f(u)+\lambda v, \quad x\in \R^N,\\
	-\Delta v +v =
	g(v)+\lambda u, \quad x\in \R^N.\end{cases}\ee
System \eqref{1-4} appears in several physical situations such as in nonlinear optics, in double Bose-Einstein condensates and in plasma physics. It has been extensively studied by many authors in the past few years. In \cite{CZ,CZ2}, Chen and Zou proved the existence of positive ground states and bound states of the coupled system \eqref{1-4} for $\lambda\in(0,1)$. More importantly, they gave more precise descriptions of the limit behavior and energy estimates of the bound states as $\lambda$ changes. In the case of $N\leq 3$, $f(s) = g(s) = s^3$, Ambrosetti, Colorado and Ruiz \cite{AmCo} proved that \eqref{1-4} has multi-bump solitons for $\lambda > 0$ small enough. When $f(u)$ and $g(v)$ are replaced by $f(x, u) = (1+c(x))|u|^{p-1}u$ and $g(x, v) = (1+d(x))|v|^{p-1}v$ respectively, system \eqref{1-4} has been studied by Ambrosetti \cite{Am} with dimension $N = 1$ and Ambrosetti, Cerami and Ruiz \cite{AmCe} with dimension $N\geq 2$. In the case of $N= 3$, $f(s) = g(s) = s^3$, Lin and Peng \cite{LinP} studied the segregated vector solutions of \eqref{1-4} as well as a 3-core coupler with circular symmetry, and by a construct argument, many positive vecrtor solutions were obtained. In \cite{LP},   L\"{u} and Peng considered a class of systems of two coupled nonlinear fractional Laplacian equations and established some results about the existence of positive vector solutions and vector ground state solutions, as well as the asymptotic behavior of these solutions as the coupling parameter tends to zero. In particular, Chen and Zou \cite{CZ1} studied the following system with one critical exponent
\be\label{1-4-0}
\begin{cases}-\Delta u +\mu u =
|u|^{q-2}u+\lambda v, \quad x\in \R^N,\\
-\Delta v +\nu v =
|v|^{2^*-2}v+\lambda u, \quad x\in \R^N,
\end{cases}
\ee
where $\mu,\nu>0$ and $0<\la<\sqrt{\mu\nu}$. They showed that system \eqref{1-4-0} has a positive ground state solution for some $\mu,\nu>0$ and $\la>0$. Moreover, if $q=2^*$, then \eqref{1-4-0} has no nontrivial solutions for $\mu,\nu>0$ and $0<\la<\sqrt{\mu\nu}$. 

\vskip0.12in
For the following system with critical exponent which are linearly coupled 
\be\label{1-4-1}
\begin{cases}-\Delta u +\mu_1u =
	|u|^{2^*-2}u+\lambda v, \quad x\in \Om,\\
	-\Delta v +\mu_2v =
	|v|^{2^*-2}v+\lambda u, \quad x\in \Om,\\
	u=v=0, \quad \mbox{on}\;\partial \Om,
\end{cases}
\ee
where $\Om$ is a smooth bounded domain in $\RN$, $N\geq 3$, $\mu_1,\mu_2>-\la_1(\Om)$, $\la_1(\Om)$ is the first eigenvalue of $(-\Delta,H_0^1(\Om))$, $\lambda\in \R$ is a coupling parameter, Peng-Shuai-Wang \cite{PSW} proved system \eqref{1-3} has a positive ground state solution for some $\la>0$ and a positive higher energy solution when $|\la|$ is small. Moreover, they analyzed the asymptotic behaviors of the positive ground state and higher energy solutions as $\la\rightarrow 0$.

For the case that $m = 1$ and $1 < p < N$ in \eqref{1-3}, the existence of a $C^1$ nonnegative solution for
\be\label{1-5}
-\Delta_pu+|u|^{p-2}u =f(u),~x\in \R^N,
\ee
has also been proved in \cite{FG} by Ferrero and Gazzola under general assumptions on $f$. In \cite{dM},  J. M. do \'{O} and E. S. Medeiros proved the existence of least energy solutions for \eqref{1-5} and established some properties of the solutions when $1<p\leq N$.

\vskip0.12in

In our paper, we assume that $f,g \in C(\R,\R)$ and are odd.

\begin{itemize}
	\item[{\bf (F1)}]  $\lim\limits_{s\rightarrow 0}\frac{f(s)}{|s|^{p-1}}=0$, $\lim\limits_{s\rightarrow 0}\frac{g(s)}{|s|^{p-1}}=0$;
	
	\item[{\bf (F2)}]
	 $\lim\limits_{|s|\rightarrow +\infty}\frac{f(s)}{|s|^{p^*-1}}=0$, $\lim\limits_{|s|\rightarrow +\infty}\frac{g(s)}{|s|^{p^*-1}}=0$;
	
	\item[{\bf (F3)}]  there exist $s_0,s_1>0$ such that $F(s_0)>\dfrac{s_0^p}{p},G(s_1)>\dfrac{s_1^p}{p}$.
\end{itemize}

The main results of the current paper are the followings.

\bt\label{th1} Suppose $N\geq 3$, $\frac{2N}{N+2}<p< 2$ and $f,g$ satisfy assumptions $(F1)-(F3)$. Then there exists $\la_0>0$ such that for $\la \in (0,\la_0)$, \eqref{1-1} has a radial solution $(u_\la,v_\la)$ with $u_\la,v_\la\in L^\iy(\RN)\cap C_{loc}^{1,\beta}$ for some $\beta\in (0,1)$. Furthermore, let $\la_n \in (0,\la_0),\;n\in \mathbb{N},$ be a sequence with $\la_n\rightarrow 0$ as $n\rightarrow \infty$. Then, passing to a subsequence, $(u_{\la_n},v_{\la_n})\rightarrow (U,V)$ strongly in $W^{1,p}(\R^N)\times W^{1,p}(\R^N)$ as $n\rightarrow \infty$, where $U$ is a positive radial ground state of \eqref{1-5} and respectively $V$ is a positive radial ground state of
\be\label{1-6}
-\Delta_pv+|v|^{p-2}v =g(v),~x\in \R^N.
\ee
\et

In the next theorem we show that we can obtain another positive vector solution for the system \eqref{1-1} which is different from the solutions obtained in Theorem \ref{th1}.

\bt\label{th2} Under the assumptions of Theorem \ref{th1}, for any $\la>0$, \eqref{1-1} has a positive radial ground state $(u_\la,v_\la)$. Furthermore, let $\la_n \in (0,1),\;n\in \mathbb{N},$ be a sequence with $\la_n\rightarrow 0$ as $n\rightarrow \infty$. Then, passing to a subsequence, $(u_{\la_n},v_{\la_n})\rightarrow (\hat{u},\hat{v})$ strongly in $W^{1,p}(\R^N)\times W^{1,p}(\R^N)$ as $n\rightarrow \infty$, and one of the following conclusions holds:
\begin{itemize}
\item[(i)]  $\hat{u}\equiv 0$ and $\hat{v}$ is a positive radial ground state of \eqref{1-6};

\item[(ii)]  $\hat{v}\equiv 0$ and $\hat{u}$ is a positive radial ground state of \eqref{1-5}.
\end{itemize}
\et

\br\label{1}
In the scalar case for $p=2$, assumptions $(F1)-(F3)$ are called Berestycki-Lions conditions, which were introduced by Berestycki and Lions \cite{BeL} to get a ground state solution for
\be\label{1-5-1}
-\Delta u+u=f(u),\;\;u\in H^1(\RN).
\ee
They showed that assumptions $(F1)-(F3)$ are almost optimal for the existence of ground states of \eqref{1-5-1} by Poho\v{z}aev identity. In our case, by the corresponding Pucci-Serrin identity \cite{DMS} for \eqref{1-5}, 
assumptions $(F1)-(F3)$ in Theorem \ref{th1} and \ref{th2} are almost optimal.
\er

To prove Theorem \ref{th1}, we use the idea from \cite{CZ2}, and we define a special mountain-pass value $c_\lambda$, where all paths are required to be bounded in $W^{1,p}(\R^N)\times W^{1,p}(\R^N)$ by the same constant which is independent of $\lambda$.

Since there is no Ambrosetti-Rabinowitz condition and the nonlinear terms $f, g$ are not homogeneous, the usual Nehari manifold method is not suitable in our case. To prove Theorems \ref{th2}, here we will adopt a minimizing argument. More precisely, let
$$P(u,v):=\int_{\R^N}\big(|\nabla u|^p+|\nabla v|^p\big)dx-p^*\int_{\R^N}\big(F(u)+G(v)+\la uv-\frac{1}{p}|u|^p-\frac{1}{p}|v|^p\big)dx,$$
and let
$$\mathcal{M}_\la:=\Big\{(u,v)\in W^{1,p}(\R^N)\times W^{1,p}\setminus \{(0,0)\}:P(u,v)=0\Big\}.$$
We give some notations here. Throughout this paper, we denote the norm of $L^p(\Om)$ by $\|u\|_p =(\int_{\Om}|u|^p\,dx)^{\frac{1}{p}}$ and positive constants (possibly different in different places) by $C, C_0, C_1, \cdots$.
Denote the norm of $W^{1,p}(\R^N)$ by
 $$\|u\|:=\big(\int_{\R^N} (|\nabla u|^p + |u|^p)\,dx\big)^{1/p},$$
 and define $W:=W^{1,p}(\R^N)\times W^{1,p}(\R^N)$ with norm $\|(u, v)\|^p:=\|u\|^p+\|v\|^p$.
Denote by $W^{1,p}_r(\R^N)$ the subspace of $W^{1,p}(\R^N)$ formed by the radially symmetric functions and define $W_r:=W^{1,p}_r(\R^N)\times W^{1,p}_r(\R^N)$.

The paper is organized as follows. In section 2,  we prove Theorems \ref{th1}. In section 3, we prove Theorem \ref{th2}.

\vskip0.3in

\s{Proof of Theorem \ref{th1}}
\renewcommand{\theequation}{2.\arabic{equation}}

In this section, we consider the functional $I_\la$ restricted to $W_r$. By Palais's Symmetric Criticality Principle in \cite{Pa}, any critical points of $I_\la: W_r\rightarrow \R$ are radially symmetric solutions of \eqref{1-1}. We assume without loss of generality that
\be\label{1-7}
f(s)\equiv g(s)\equiv 0~~~\hbox{for\;all}\; s\leq 0.
\ee
The energy functionals of \eqref{1-5} and \eqref{1-6} are given by
$$J_1(u)=\frac{1}{p}\int_{\R^N}\big(|\nabla u|^p+|u|^p\big)dx-\int_{\R^N}F(u)dx,\;\;u\in W^{1,p}(\R^N),$$
$$J_2(v)=\frac{1}{p}\int_{\R^N}\big(|\nabla v|^p+|v|^p\big)dx-\int_{\R^N}G(v)dx,\;\;v\in W^{1,p}(\R^N).$$
Under assumptions $(F1)-(F3)$, the authors in \cite{dM} proved that \eqref{1-5} (resp.
\eqref{1-6}) has a ground state solution, and each solution $U$ of \eqref{1-5} (resp. $V$ of
\eqref{1-6}) satisfies Poho\v{z}aev-Pucci-Serrin identity:
\be\label{2-1}J_3(U)=:(N-p)\int_{\R^N}|\nabla U|^pdx+N\int_{\R^N}|U|^pdx-Np\int_{\R^N}F(U)dx=0,\ee
\be\label{2-2}J_4(V)=:(N-p)\int_{\R^N}|\nabla V|^pdx+N\int_{\R^N}|V|^pdx-Np\int_{\R^N}G(V)dx=0.\ee
Byeon, Jeanjean and Maris \cite{BJM} proved that ground state solutions of \eqref{1-5} (resp.\eqref{1-6}) must be radial up to a translation. By \eqref{1-7} and the strong maximum principle in \cite[Theorem 5]{V}, any nontrivial solutions
of \eqref{1-5} (resp.\eqref{1-6}) must be positive. Define
$$S_1:=\{U\in W^{1,p}_r(\R^N): U~ \hbox{is a positive ground state of}~ \eqref{1-5}\},$$
$$S_2:=\{V\in W^{1,p}_r(\R^N): V~\hbox{is a positive ground state of}~ \eqref{1-6}\},$$
$$X:=S_1\times S_2.$$
Take fixed $U_0 \in S_1$ and $V_0 \in S_2$, and denote the least energy of \eqref{1-5} and \eqref{1-6} respectively by
$$M_1:=J_1(U_0)~~\hbox{and}~~M_2:=J_2(V_0).$$
Then $M_1>0,\;M_2>0$. Moreover, \cite[Lemma 2.4]{dM} says that
\be\label{2-3}
M_1=\inf_{\substack{u\in W^{1,p}(\R^N)\setminus \{0\} \\ J_3(u)=0}}J_1(u),\;\;M_2=\inf_{\substack{v\in W^{1,p}(\R^N)\setminus \{0\} \\ J_4(v)=0}}J_2(v).
\ee
Without loss of generality, we assume that $M_1\leq M_2$.

\vskip0.2in
\bl\label{lemma2-1} Under the assumptions of Theorem \ref{th1}, we have
\begin{itemize}
\item[(1)]  for any $(u,v)\in X$, $u,v\in C^{1,\al}_{loc}(\RN)\cap L^\iy(\RN)$ for some $\al\in (0,1)$;

\item[(2)]  X is compact in $W$, and there exist constants $C_2 > C_1 > 0$ such
that
$C_1 \le \|(u,v)\| \leq C_2,\;\;\forall (u,v)\in X$.
\end{itemize}
\el

\noindent {\bf Proof. }  (1) The boundedness of $L^\iy$ norm of $u$ and $v$ is similar to \cite{LY} and we omit the proof. By \cite{T}, we deduce that $u,\;v\in C^{1,\al}_{loc}(\RN)$ for some $0<\al<1$.

(2) For $u\in S_1$, by \eqref{2-1}, we have
\be\label{2-4}
\begin{split}
M_1=J_1(u)=&\frac{1}{p}\int_{\R^N}\big(|\nabla u|^p+|u|^p\big)dx-\int_{\R^N}F(u)dx\\
=&(\frac{1}{p}-\frac{1}{p^*})\int_{\R^N}\big|\nabla u|^pdx\\
=&\frac{1}{N}\int_{\R^N}\big|\nabla u|^pdx.
\end{split}
\ee
Then we get that $\{\|\nabla u\|_p:u\in S_1\}$ is bounded, which implies $\{\|u\|_{p^*}:u\in S_1\}$ is also bounded. By $(F1)-(F2)$, there exists $C>0$ such that $|F(s)|\leq \frac{1}{2p}|s|^p+C|s|^{p^*}$. Then
\be\label{2-4-1}
\begin{split}
M_1=&\frac{1}{p}\int_{\R^N}\big(|\nabla u|^p+|u|^p\big)dx-\int_{\R^N}F(u)dx\\
\geq &\frac{1}{p}\int_{\R^N}\big|\nabla u|^pdx+\frac{1}{2p}\int_{\R^N}|u|^pdx-C\int_{\R^N}|u|^{p^*}dx,
\end{split}
\ee
which implies that $\{\|u\|_p:u\in S_1\}$ is bounded as well. Then we have that $S_1$ is bounded in $W^{1,p}(\R^N)$. For any sequence $\{u_n\} \subset S_1$, we can assume $u_n\rightharpoonup u_0$ in $W^{1,p}(\R^N)$ and $u_n\rightarrow u_0$ in $L_{loc}^q(\R^N)$ where $p\leq q<p^*$ up to a subsequence. We need to show that $u_n\rightarrow u_0$ in $W^{1,p}(\R^N)$. Indeed, since $J_1'(u_n)=0$ in $W^{-1,p'}(\R^N)$, for any $\varphi\in C_c^\iy(\RN)$, we have
$$\int_{\R^N}|\nabla u_n|^{p-2}\nabla u_n\nabla (u_0\varphi)+|u_n|^{p-2}u_n(u_0\varphi)=\int_{\R^N}f(u_n)(u_0\varphi),$$
$$\int_{\R^N}|\nabla u_n|^{p-2}\nabla u_n\nabla (u_n\varphi)+|u_n|^{p-2}u_n(u_n\varphi)=\int_{\R^N}f(u_n)(u_n\varphi).$$
Then we have
\be\label{2-4-2}
\begin{split}
&\int_{\R^N}|\nabla u_n|^{p-2}\nabla u_n(\nabla u_n-\nabla u_0)\varphi\\
=&\int_{\R^N}|\nabla u_n|^{p-2}\nabla u_n[\nabla (u_n\varphi)-\nabla (u_0\varphi)-(u_n-u_0)\nabla \varphi]\\
=&\int_{\R^N}f(u_n)(u_n-u_0)\varphi-\int_{\R^N}|u_n|^{p-2}u_n(u_n-u_0)\varphi-\int_{\R^N}|\nabla u_n|^{p-2}\nabla u_n(u_n-u_0)\nabla \varphi.
\end{split}
\ee
Since $u_n\rightharpoonup u_0$ in $W^{1,p}(\R^N)$ and $u_n\rightarrow u_0$ in $L_{loc}^q(\R^N)$ where $p\leq q<p^*$, by $(F1)-(F2)$, for any $\varepsilon>0$, there exists $C_\varepsilon>0$ such that $|f(s)|\leq \varepsilon |s|^{p-1}+\varepsilon |s|^{p^*-1}+C_\varepsilon |s|$. Then it's easy to check that
$$\int_{\R^N}f(u_n)(u_n-u_0)\varphi\rightarrow 0,$$
$$\int_{\R^N}|u_n|^{p-2}u_n(u_n-u_0)\varphi\rightarrow 0,$$
$$\int_{\R^N}|\nabla u_n|^{p-2}\nabla u_n(u_n-u_0)\nabla \varphi\rightarrow 0,$$
and
$$\int_{\R^N}|\nabla u_0|^{p-2}\nabla u_0(\nabla u_n-\nabla u_0)\varphi\rightarrow 0.$$
Therefore, we have
$$\int_{\R^N}\big(|\nabla u_n|^{p-2}\nabla u_n-|\nabla u_0|^{p-2}\nabla u_0,\nabla u_n-\nabla u_0\big)\varphi\rightarrow 0.$$
Using a well known inequality found in \cite{AY} \cite[Lemma A.0.5]{Pe}, we know that
$$\big(|\xi|^{p-2}\xi-|\eta|^{p-2}\eta,\xi-\eta\big)\geq\begin{cases}d_1|\xi-\eta|^p,\;\;&\hbox{if}~~p\geq2,\\
d_2(|\xi|+|\eta|)^{p-2}|\xi-\eta|^2,\;&\hbox{if}~~p\in (1,2),\end{cases}$$
where $d_1,d_2$ are positive constants. For $p\in (1,2)$, it follows that
\be\label{2-4-3}
\begin{split}
&\big(\int_{\R^N}|\nabla u_n-\nabla u_0|^p\varphi\big)^\frac{2}{p}\\
=&\big(\int_{\R^N}\frac{|\nabla u_n-\nabla u_0|^p}{(|\nabla u_n|+|\nabla u_0|)^\frac{p(2-p)}{2}}\varphi^\frac{p}{2}(|\nabla u_n|+|\nabla u_0|)^\frac{p(2-p)}{2}\varphi^\frac{2-p}{2}\big)^\frac{2}{p}\\
\leq&\big(\int_{\R^N}\frac{|\nabla u_n-\nabla u_0|^2}{(|\nabla u_n|+|\nabla u_0|)^{2-p}}\varphi\big)\big(\int_{\R^N}(|\nabla u_n|+|\nabla u_0|)^p\varphi\big)^\frac{2-p}{p}\\
\leq&C\int_{\R^N}\big(|\nabla u_n|^{p-2}\nabla u_n-|\nabla u_0|^{p-2}\nabla u_0,\nabla u_n-\nabla u_0\big)\varphi\rightarrow 0.
\end{split}
\ee
Similarly, we can prove the same local convergence property for the case $p\geq 2$. Then we can deduce that
$$\nabla u_n\rightarrow \nabla u_0~~a.e.~~x\in \RN.$$

Now consider the following minimizing problem
$$T_\mu:=\inf\limits_{u\in W^{1,p}(\R^N)}\Big\{\frac{1}{p}\int_{\R^N}|\nabla u|^pdx:\int_{\R^N}F(u)dx-\frac{1}{p}\int_{\R^N}|u|^pdx=\mu\Big\},$$
where $\mu>0$, and it has been shown that there exists a minimizer for $T_1$ in \cite[Theorem 1.4]{dM}. Then by \cite[Lemma 1]{BJM}, we know that $u_n$ is a minimizer for $T_{\mu_0}$, where $\mu_0 =(\frac{N-p}{N}T_1)^{N/p}$. Then $\{u_n\}$ is a minimizing sequence for $T_{\mu_0}$ and $u_n$ is positive and radially symmetric. As in \cite{BeL}, we know that $u_0$ is a minimizer for $T_{\mu_0}$. Then we have
$$\int_{\R^N}F(u_0)-\frac{1}{p}|u_0|^p=\int_{\R^N}F(u_n)-\frac{1}{p}|u_n|^p=\mu_0,$$
$$\frac{1}{p}\int_{\R^N}|\nabla u_n|^p=\frac{1}{p}\int_{\R^N}|\nabla u_0|^p= T_{\mu_0}.$$
Since $\nabla u_n\rightarrow \nabla u_0~~a.e.~~x\in \RN$, by Brezis-Lieb lemma, we have that $\nabla u_n\rightarrow \nabla u_0$ in $L^p(\RN)$. Since $W_r^{1,p}(\R^N)\hookrightarrow L^q(\RN)$ is compact for $p<q<p^*$, by $(F1)-(F2)$, we have $\int_{\R^N}F(u_n)\rightarrow \int_{\R^N}F(u_0)$. Then by Brezis-Lieb lemma again, we deduce that $u_n\rightarrow u_0$ in $L^p(\RN)$. Therefore we get that $u_n\rightarrow u_0$ in $W^{1,p}(\R^N)$, which implies that $S_1$ is compact in $W^{1,p}(\R^N)$. Similarly, $S_2$ is compact in $W^{1,p}(\R^N)$.
\qed

\vskip0.2in
For $t, s > 0$, we define $U_{0,t}(x) := U_0( \frac{x}{t})$ and $V_{0,s}(x) := V_0(\frac{x}{s})$. Then by \eqref{2-1} we have
\be\label{2-5}
\begin{split}
J_1(U_{0,t})=&\frac{t^{N-p}}{p}\int_{\R^N}|\nabla U_0|^pdx+\frac{t^N}{p}\int_{\R^N}|U_0|^pdx-t^N\int_{\R^N}F(U_0)dx\\
=&\big(\frac{t^{N-p}}{p}-\dfrac{(N-p)t^N}{Np}\big)\int_{\R^N}|\nabla U_0|^pdx.
\end{split}
\ee
Note that
\be\label{2-6}
J_1(U_{0,1})=\max\limits_{t>0}J_1(U_{0,t})=J_1(U_0)=M_1.
\ee
It is easily seen that there exists $0 < t_0 < 1 < t_1$ such that
\be\label{2-7}
J_1(U_{0,t})\leq \dfrac{1}{4}M_1\;\;\hbox{for}\;t\in (0,t_0]\cup [t_1,\infty).
\ee
Similarly, there exists $0 < s_0 < 1 < s_1$ such that
\be\label{2-8}
J_2(V_{0,s})\leq \dfrac{1}{4}M_1\;\;\hbox{for}\;s\in (0,s_0]\cup [s_1,\infty).
\ee
Define
$$\tilde{\ga}_1(t):=U_{0,t}\;\;\hbox{for}\;0<t\leq t_1,\;\;\tilde{\ga}_1(0):=0;$$
$$\tilde{\ga}_2(s):=V_{0,s}\;\;\hbox{for}\;0<s\leq s_1,\;\;\tilde{\ga}_2(0):=0;$$
$$\tilde{\ga}(t,s):=(\tilde{\ga}_1(t),\tilde{\ga}_2(s)).$$
Then there exists a constant $\mathcal{C} > 0$ such that
\be\label{2-9}
\max\limits_{(t,s)\in Q}\|\tilde{\ga}(t,s)\|=\max\limits_{(t,s)\in Q}\int_{\R^N}t^{N-p}|\nabla U_0|^p+t^N|U_0|^p+\int_{\R^N}s^{N-p}|\nabla V_0|^p+s^N|V_0|^p \leq \mathcal{C},
\ee
where $Q:=[0,t_1]\times [0,s_1]$. Recalling $C_2$ in Lemma \ref{lemma2-1}, we define
$$c_\la:=\inf\limits_{\ga\in \Gamma}\max\limits_{(t,s)\in Q}I_\la(\ga(t,s)),\;\;d_\la:=\max\limits_{(t,s)\in Q}I_\la(\tilde{\ga}(t,s)),$$
where
\begin{align}\label{2-10}
\Gamma:=\big\{&\ga\in C(Q,W_r):\max\limits_{(t,s)\in Q}\|\ga(t,s)\| \leq 2C_2+\mathcal{C},\nonumber\\
&\ga(t,s)=\tilde{\ga}(t,s)\;\;\hbox{for}\;(t,s)\in Q\setminus (t_0,t_1)\times (s_0,s_1)\big\}.
\end{align}

\vskip0.2in
\bl\label{lemma2-2} $\lim\limits_{\la \rightarrow 0}c_\la=\lim\limits_{\la \rightarrow 0}d_\la=c_0=d_0=M_1+M_2$.
\el

\noindent {\bf Proof. } Since $\la\geq 0$, by the definition of $I_\la$ in \eqref{1-2}, we have
$I_\la(\tilde{\ga}(t,s))\leq I_0(\tilde{\ga}(t,s))$, and so
\be\label{2-11}
\begin{split}
d_\la&\leq d_0=\max\limits_{(t,s)\in Q}I_0(\tilde{\ga}(t,s))=\max\limits_{t\in [0,t_1]}J_1(\tilde{\ga}_1(t))+\max\limits_{s\in [0,s_1]}J_2(\tilde{\ga}_2(s))\\
&=J_1(\tilde{\ga}_1(1))+J_2(\tilde{\ga}_2(1))=J_1(U_0)+J_2(V_0)=M_1+M_2.
\end{split}
\ee
Since $\tilde{\ga}\in \Gamma$, we have $c_\la\leq d_\la$, and then
\be\label{2-12}
\limsup\limits_{\la\rightarrow 0}c_\la \leq \limsup\limits_{\la\rightarrow 0}d_\la\leq d_0,\;\;c_0\leq d_0.
\ee
On the other hand, for any $\ga(t,s)=(\ga_1(t,s),\ga_2(t,s))\in \Gamma$, we define $\Upsilon(\ga):[t_0,t_1]\times [s_0,s_1]\rightarrow \R^2$ by
$$\Upsilon(\ga)(t,s):=\big(J_5(\ga_1(t,s))-1,J_6(\ga_2(t,s))-1\big),$$
where $J_5,J_6:W^{1,p}(\R^N)\rightarrow\R$ are defined by
$$J_5(u):=\begin{cases}\frac{Np\int_{\R^N}F(u)dx}{\int_{\R^N}(N-p)|\nabla u|^p+N|u|^pdx},\;\;&u\neq 0,\\
0,\;&u=0,\end{cases},$$
$$J_6(u):=\begin{cases}\frac{Np\int_{\R^N}G(u)dx}{\int_{\R^N}(N-p)|\nabla u|^pdx+N|u|^pdx},\;\;&u\neq 0,\\
0,\;&u=0.\end{cases}$$
By $(F1)-(F2)$ and the Sobolev inequality, it is easy to prove that
$J_5, J_6$ are continuous. Similarly as in \eqref{2-5} it is easily seen that
\be\label{2-12-1}
\Upsilon(\tilde{\ga})(t,s)=\Big(\dfrac{p^*t^p\int_{\R^N}F(U_0)dx}{\int_{\R^N}|\nabla U_0|^p+p^*/pt^p|U_0|^p}-1,
\dfrac{p^*s^p\int_{\R^N}G(V_0)dx}{\int_{\R^N}|\nabla V_0|^p+p^*/ps^p|V_0|^p}-1\Big).
\ee
Recalling \eqref{2-1} and \eqref{2-2}, we have $\Upsilon(\tilde{\ga})(1,1)=(0,0)$. By a direct computation, one gets that $\deg(\Upsilon(\tilde{\ga}), [t_0, t_1]\times [s_0, s_1], (0, 0)) = 1$. By \eqref{2-10},
we see that for any $(t, s) \in \partial([t_0, t_1]\times [s_0, s_1]), \Upsilon(\ga)(t,s)= \Upsilon(\tilde{\ga})(t,s)\neq (0, 0)$.
Therefore, $\deg(\Upsilon(\ga)), [t_0, t_1]\times [s_0, s_1], (0, 0))$ is well defined and
$$\deg(\Upsilon(\ga)), [t_0, t_1]\times [s_0, s_1], (0, 0))=\deg(\Upsilon(\tilde{\ga}), [t_0, t_1]\times [s_0, s_1], (0, 0))=1.$$
Then there exists $(t_2, s_2) \in [t_0, t_1]\times [s_0, s_1]$ such that $\Upsilon(\ga)(t_2, s_2) = (0, 0)$, that is,
$J_5(\ga_1(t_2, s_2)) = J_6(\ga_2(t_2, s_2)) = 1$. This implies $J_3(\ga_1(t_2, s_2)) = J_4(\ga_2(t_2, s_2)) = 0$
and $\ga_i(t_2, s_2)\neq 0$ for $i = 1, 2$. Combining these with \eqref{2-3}, we have
\be\label{2-13}
\begin{split}
\max\limits_{(t,s)\in Q}I_0(\ga(t,s))&\geq I_0(\ga(t_2, s_2))=J_1(\ga_1(t_2, s_2))+J_2(\ga_2(t_2, s_2))\\
&\geq M_1+M_2=d_0.
\end{split}\ee
Therefore, $c_0\geq d_0$. By \eqref{2-12}, we have $c_0=d_0$.

Finally, assume by contradiction that $\liminf\limits_{\lambda\rightarrow 0}c_\lambda<d_0$. Then there exists $\varepsilon> 0, \lambda_n\rightarrow 0$ and $\ga_n=(\ga_{n,1}, \ga_{n,2})\in \Gamma$ such that
$$\max\limits_{(t,s)\in Q}I_{\lambda_n}(\ga_n(t,s))\leq d_0-2\varepsilon.$$
By the definition of $\Gamma $ in \eqref{2-10} and H\"{o}lder's inequality, there exist $C>0$ and $n_0$ large enough such that
$$\max\limits_{(t,s)\in Q}\lambda_n\big|\int_{\R^N}\ga_{n,1}(t,s)\ga_{n,2}(t,s)dx\big|\leq C\lambda_n\leq \varepsilon,\;\forall n\geq n_0,$$
and then
$$\max\limits_{(t,s)\in Q}I_0(\ga_n(t,s))\leq \max\limits_{(t,s)\in Q}I_{\lambda_n}(\ga_n(t,s))+\varepsilon\leq d_0-\varepsilon,\;\forall n\geq n_0,$$
a contradiction with \eqref{2-13}. Therefore, $\liminf_{\lambda\rightarrow 0}c_\lambda\geq d_0$. Combining this with \eqref{2-12}, we complete the proof.
\qed
\vskip0.2in
Recalling that $X=S_1\times S_2$, we define
$$X^\delta :=\{(u,v)\in W_r:\hbox{dist}((u,v),X)\leq \delta\},\;\;I_\la^d:=\{(u,v)\in W_r:I_\la(u,v)\leq d\}.$$

\bl\label{lemma2-3} Let $C_1$ be in Lemma \ref{lemma2-1}. For a small $\delta\in (0,C_1/2)$, there exists constants $0<\sigma<1$ and $\la_1>0$ such that $\|I_\la'(u,v)\|\geq \sigma$ for any $(u,v)\in I_\lambda^{d_\lambda}\cap (X^\delta\setminus X^{\delta/2})$ and any $\lambda \in (0,\lambda_1)$.
\el

\noindent {\bf Proof. } Assume by contradiction that there exist $\lambda_n\rightarrow 0$ and $(u_n,v_n)\in I_{\lambda_n}^{d_{\lambda_n}}\cap (X^\delta\setminus X^{\delta/2})$ such that $\|I_{\la_n}'(u_n,v_n)\|\rightarrow 0$. By Lemma \ref{lemma2-1}, $\{(u_n,v_n),n\geq 1\}$ are uniformly bounded in $W_r$. Recall that the Sobolev embedding $W^{1,p}_r(\R^N)\hookrightarrow L^q(\R^N)$ is compact for any $q\in (p,p^*)$. Up to a subsequence, we may assume that $(u_n,v_n)\rightharpoonup (U,V)$ weakly in $W_r$ and strongly in  $L^{q_1}(\R^N)\times L^{q_2}(\R^N),\;q_1,q_2\in (p,p^*)$. As in the proof of Lemma \ref{lemma2-1}, we get that $\nabla u_n\rightarrow \nabla U$ a.e. $x\in \RN$. Let $\psi\in C_c^\iy(\RN)$ with $\psi(x)=1, |x|\leq 1$ and $\psi(x)=0, |x|\geq 2$. As in \cite[Theorem A.4]{Willem.1996}, we define $h_1(u)=\psi(u)f(u), h_2(u)=(1-\psi(u))f(u)$, then by $(F1)-(F2)$, we have $|h_1(u)|\leq C|u|^{p-1}, |h_2(u)|\leq C|u|^{p^*-1}$. Therefore, we get that $\{|u_n|^{p-2}u_n\}, \{|\nabla u_n|^{p-2}\nabla u_n\}, \{h_1(u_n)\}$ are bounded in $L^\frac{p}{p-1}(\RN)$ and $\{h_2(u_n)\}$ is bounded in $L^\frac{p^*}{p^*-1}(\RN)$. By \cite{BL}, we have that $|u_n|^{p-2}u_n\rightharpoonup |U|^{p-2}U, |\nabla u_n|^{p-2}\nabla u_n\rightharpoonup |\nabla U|^{p-2}\nabla U, h_1(u_n)\rightharpoonup h_1(U)$ in $L^\frac{p}{p-1}(\RN)$ and $h_2(u_n)\rightharpoonup h_2(U)$ in $L^\frac{p^*}{p^*-1}(\RN)$. Then for any $\varphi\in W^{1,p}_r(\R^N)$, we have
$$I_{\la_n}'(u_n,v_n)(\varphi,0)\longrightarrow I_0'(U,V)(\varphi,0)=0.$$
Similarly we can get the same results for $v_n$ and then
$I_0'(U,V)=0$; that is, $U$ (resp. $V$) is a solution of \eqref{1-5} (resp. \eqref{1-6}). Moreover, since $u_n\rightarrow U$ in $L^q(\R^N),\;q\in (p,p^*)$, by $(F1)-(F2)$ again, we have
$$\lim\limits_{n\rightarrow \iy}\int_{\R^N}f(u_n)u_ndx=\int_{\R^N}f(U)Udx.$$
By $I_{\la_n}'(u_n,v_n)(u_n,0)\rightarrow 0$, we get that
$$\|u_n\|^p=\int_{\R^N}f(u_n)u_ndx+o(1)=\int_{\R^N}f(U)Udx+o(1)=\|U\|^p+o(1),$$
and so $u_n\rightarrow U$ strongly in $W^{1,p}_r(\R^N)$. Similarly, $v_n\rightarrow V$ strongly in $W^{1,p}_r(\R^N)$, and so $(U,V)\in X^\delta$, which implies that $U\not\equiv​ 0$ and $V\not\equiv​0$. By \eqref{1-7} and the strong maximum principle, we have $U,V>0$. Recalling Lemma \ref{lemma2-2} and the definition of $M_1,M_2$, we have
\begin{align}
M_1+M_2&\leq J_1(U)+J_2(V)=I_0(U,V)=\lim\limits_{n\rightarrow \iy}I_{\la_n}(u_n,v_n)\nonumber\\
&\leq \lim\limits_{n\rightarrow \iy}d_{\la_n}=M_1+M_2.\nonumber
\end{align}
This implies $J_1(U)=M_1,\;J_2(V)=M_2$, that is, $U\in S_1,\; V\in S_2$. So $(U,V)\in X$, which contradicts with $(u_n,v_n)\notin X^{\delta/2}$ for any $n\geq 1$. This completes the proof.
\qed

\vskip0.2in
From now on, we fix a small $\delta \in(0, \min{\mathcal{C}/2, C_1/2})$ and corresponding $0 <\sigma< 1$ and $\lambda_1 > 0$ such that conclusions in Lemma \ref{lemma2-3} hold.

\bl\label{lemma2-4} There exist $\la_2\in (0,\la_1]$ and $\al>0$ such that for any $\la\in(0,\la_2)$,
$$I_\la(\tilde{\ga}(t,s))\geq c_\la-\al\;\;implies~that~~\tilde{\ga}(t,s)\in  X^{\delta/2}.$$
\el

\noindent {\bf Proof. } Assume by contradiction that there exist $\la_n\rightarrow 0,\;\al_n\rightarrow 0$ and $(t_n,s_n)\in Q$ such that
\be\label{2-14}
I_{\la_n}(\tilde{\ga}(t_n,s_n))\geq c_{\la_n}-\al_n\;\;\hbox{and}\;\; \tilde{\ga}(t_n,s_n)\notin  X^{\delta/2}. \ee
Passing to a subsequence, we may assume that $(t_n,s_n)\rightarrow (\tilde{t},\tilde{s})\in Q$. Then by Lemma \ref{lemma2-2}, we have
$$I_0(\tilde{\ga}(\tilde{t},\tilde{s}))=\lim\limits_{n\rightarrow \iy}I_{\la_n}(\tilde{\ga}(t_n,s_n))\geq \lim\limits_{n\rightarrow \iy}c_{\la_n}=M_1+M_2.$$
Combining this with \eqref{2-5},\eqref{2-6} and \eqref{2-11}, it's easy to see that $ (\tilde{t},\tilde{s})=(1,1).$ Hence,
$$\lim\limits_{n\rightarrow \iy}\|\tilde{\ga}(t_n,s_n)-\tilde{\ga}(1,1)\|=0.$$
However, $\tilde{\ga}(1,1)=(U_0,V_0)\in X$, which is a contradiction to \eqref{2-14}.\qed

\vskip0.2in
Let
\be\label{2-15}
\al_0:=\min\big\{\frac{\al}{2},\frac{M_1}{4},\frac{1}{8}\delta \sigma^2\big\},
\ee
where $\delta,\sigma$ are seen in Lemma \ref{lemma2-3}. By Lemma \ref{lemma2-2}, there exists $\la_0\in (0,\la_2]$ such that
\be\label{2-16}
|c_\la-d_\la|<\al_0,\;\;|c_\la-(M_1+M_2)|<\al_0,\;\forall \la\in (0,\la_0).
\ee

\bl\label{lemma2-5} For fixed $\la\in (0,\la_0)$, there exists $\{(u_n,v_n)\} \subset X^\delta\cap I_\la^{d_\la}$ such that
$$I_\la'(u_n,v_n)\rightarrow 0~~in~~W_r~~as~~n\rightarrow \iy.$$
\el

\noindent {\bf Proof. } Fix a $\la\in (0,\la_0)$. Assume by contradiction that there exists $0<l(\la)<1$ such that $\|I_\la'(u,v)\|\geq l(\la)$ on $X^\delta\cap I_\la^{d_\la}$. Then there exists a locally Lipschitz continuous pseudo-gradient vector field $T_\la$ in $W_r$ which is defined on a neighborhood $Z_\la$ of $X^\delta\cap I_\la^{d_\la}$ (see \cite[Lemma 3.2]{St}) such that for any $(u,v)\in Z_\la$, there holds
$$\|T_\la(u,v)\|\leq 2\min\{1,\|I_\la'(u,v)\|\},$$
$$\langle I_\la'(u,v),T_\la(u,v)\rangle\geq \min\{1,\|I_\la'(u,v)\|\}\|I_\la'(u,v)\|.$$
Let $\eta_\la$ be a Lipschitz continuous function on $W_r$ such that $0\leq \eta_\la\leq 1,\;\eta_\la\equiv 1$ on $X^\delta\cap I_\la^{d_\la}$ and $\eta_\la\equiv 0$ on $W_r\setminus Z_\la$. Let $\xi_\la$ be a Lipschitz continuous function on $\R$ such that $0\leq \xi_\la\leq 1,\;\xi_\la(s)\equiv 1$ if $|s-c_\la|\leq \frac{\al}{2}$ and $\xi_\la(s)\equiv 0$ if $|s-c_\la|\geq \al$. Let
$$e_\la(u,v):=\begin{cases}-\eta_\la(u,v)\xi_\la(I_\la(u,v))T_\la(u,v),\;\;&(u,v)\in Z_\la,\\
0,\;&(u,v)\in W_r\setminus Z_\la.\end{cases}$$
It's easy to see that $e_\la$ is locally Lipschitz continuous throughout $W_r$. Moreover, since $\|T_\la(u,v)\|\leq 2$ uniformly, also $\|e_\la(u,v)\|\leq 2$ is uniformly bounded. Then there exists a global solution $\psi_\la:W_r\times [0,+\iy)\rightarrow W_r$ for the initial value problem
$$\begin{cases}-\frac{d}{d\theta}\psi_\la(u,v,\theta)=e_\la(\psi_\la(u,v,\theta)),\\
\psi_\la(u,v,0)=(u,v).\end{cases}$$
And $\psi_\la$ has the following properties:
\begin{itemize}
\item[(1)] $\psi_\la(u,v,\theta)=(u,v)$ if $\theta=0$ or $(u,v)\in W_r\setminus Z_\la$ or $|I_\la(u,v)-c_\la|\geq \al$;
\item[(2)] $\|\frac{d}{d\theta}\psi_\la(u,v,\theta)\|\leq 2$;
\item[(3)] $\frac{d}{d\theta}I_\la(\psi_\la(u,v,\theta))=\langle I_\la'(\psi_\la(u,v,\theta)), e_\la(\psi_\la(u,v,\theta))\rangle\leq 0$.
\end{itemize}
{Step 1.} For any $(t,s)\in Q$, we claim that there exists $\theta_{t,s}\in [0,+\iy)$ such that $\psi_\la(\tilde{\ga}(t,s),\theta_{t,s})\in I_\la^{c_\la-\al_0}$, where $\al_0$ is seen in \eqref{2-15}.

Assume by contradiction that there exists $(t,s)\in Q$ such that  $$I_\la(\psi_\la(\tilde{\ga}(t,s),\theta))>c_\la-\al_0,$$
for any $\theta \geq 0$. Note that $\al_0<\al$; we see from Lemma \ref{lemma2-4} that $\tilde{\ga}(t,s)\in X^{\delta/2}$. Note that $I_\la(\tilde{\ga}(t,s))\leq d_\la<c_\la+\al_0$; we see from the property (3) that
$$c_\la-\al_0<I_\la(\psi_\la(\tilde{\ga}(t,s),\theta))\leq d_\la<c_\la+\al_0,\;\forall \theta\geq 0.$$
This implies $\xi_\la(I_\la(\psi_\la(\tilde{\ga}(t,s),\theta)))\equiv 1$. If $\psi_\la(\tilde{\ga}(t,s),\theta)\in X^\delta$ for all $\theta\geq 0$, then $\eta_\la(\psi_\la(\tilde{\ga}(t,s),\theta))\equiv 1$, and $\|I_\la'(\psi_\la(\tilde{\ga}(t,s),\theta))\|\geq l(\la)$ for all $\theta>0$. Then
\begin{align*}
I_\la\Big(\psi_\la\big(\tilde{\ga}(t,s),\frac{\al}{l(\la)^2}\big)\Big)&\leq I_\la\Big(\psi_\la\big(\tilde{\ga}(t,s),0\big)\Big)+\int_{0}^{\frac{\al}{l(\la)^2}}\frac{d}{d\theta}I_\la\Big(\psi_\la\big(\tilde{\ga}(t,s),\theta\big)\Big)d\theta\\
&\leq c_\la+\al_0-\int_{0}^{\frac{\al}{l(\la)^2}}l(\la)^2d\theta\\
&=c_\la+\al_0-\al\\
&\leq c_\la-\al_0,
\end{align*}
a contradiction. Thus, there exists $\theta_{t,s}>0$ such that $\psi_\la(\tilde{\ga}(t,s),\theta_{t,s})\notin X^\delta$. Note that $\tilde{\ga}(t,s) \in X^{\delta/2}$; there exist $0<\theta_{t,s}^1<\theta_{t,s}^2\leq\theta_{t,s}$ such that $\psi_\la(\tilde{\ga}(t,s),\theta_{t,s}^1)\in \partial X^{\delta/2},\;\psi_\la(\tilde{\ga}(t,s),\theta_{t,s}^2)\in \partial X^\delta$ and $\psi_\la(\tilde{\ga}(t,s),\theta)\in X^\delta\setminus X^{\delta/2}$ for all $\theta\in (\theta_{t,s}^1,\theta_{t,s}^2)$. Then be Lemma \ref{lemma2-3}, we have
$\|I_\la'(\psi_\la(\tilde{\ga}(t,s),\theta))\|\geq \sigma$ for all $\theta\in (\theta_{t,s}^1,\theta_{t,s}^2)$. Then using property (2), we have
$$\delta/2\leq \|\psi_\la(\tilde{\ga}(t,s),\theta_{t,s}^2)-\psi_\la(\tilde{\ga}(t,s),\theta_{t,s}^1)\|\leq 2|\theta_{t,s}^2-\theta_{t,s}^1|,$$
that is, $\theta_{t,s}^2-\theta_{t,s}^1\geq \delta/4$. This implies
\begin{align*}
I_\la\Big(\psi_\la\big(\tilde{\ga}(t,s),\theta_{t,s}^2\big)\Big)&\leq I_\la\Big(\psi_\la\big(\tilde{\ga}(t,s),\theta_{t,s}^1\big)\Big)+\int_{\theta_{t,s}^1}^{\theta_{t,s}^2}\frac{d}{d\theta}I_\la\Big(\psi_\la\big(\tilde{\ga}(t,s),\theta\big)\Big)d\theta\\
&< c_\la+\al_0-\sigma^2(\theta_{t,s}^2-\theta_{t,s}^1)\\
&\leq c_\la+\al_0-\frac{1}{4}\delta\sigma^2\\
&\leq c_\la-\al_0,
\end{align*}
which is a contradiction.

By {Step 1.}, we can define $T(t,s):=\inf\{\theta\geq 0:I_\la(\psi_\la(\tilde{\ga}(t,s),\theta))\leq c_\la-\al_0\}$ and let $\ga(t,s):=\psi_\la(\tilde{\ga}(t,s), T(t,s))$. Then $I_\la(\ga(t,s))\leq c_\la-\al_0$ for all $(t,s)\in Q$.

{Step 2.} We shall prove that $\ga(t,s)\in \Gamma$.

For any $(t,s)\in Q\setminus (t_0,t_1)\times (s_0,s_1)$, by \eqref{2-4}, \eqref{2-6}-\eqref{2-8} and \eqref{2-15}-\eqref{2-16}, we have
\begin{align*}
I_\la(\tilde{\ga}(t,s))&\leq I_0(\tilde{\ga}(t,s))=J_1(\tilde{\ga}_1(t))+J_2(\tilde{\ga}_2(s))\\
&\leq \frac{M_1}{4}+M_2\leq M_1+M_2-3\al_0<c_\la-\al_0,
\end{align*}
which implies that $T(t,s)=0$ and so $\ga(t,s)=\tilde{\ga}(t,s)$.

From the definition of $\Gamma$ in \eqref{2-10}, it suffices to prove that
$\|\ga(t,s)\|\leq 2C_2+\mathcal{C}$ for all $(t, s) \in Q$ and $T(t, s)$ is continuous with respect to $(t, s)$.

For any $(t,s)\in Q$, if $I_\la(\tilde{\ga}(t,s))\leq c_\la-\al_0$, we have $T(t,s)=0$ and so $\ga(t,s)=\tilde{\ga}(t,s)$, and by \eqref{2-9}, we see that $\|\ga(t,s)\|\leq \mathcal{C}< 2C_2+\mathcal{C}$.

If $I_\la(\tilde{\ga}(t,s))> c_\la-\al_0$, then $\tilde{\ga}(t,s)\in X^{\delta/2}$ and
$$c_\la-\al_0<I_\la(\psi_\la(\tilde{\ga}(t,s),\theta))\leq d_\la<c_\la+\al_0,\;\;\forall \theta \in [0,T(t, s)).$$
This implies $\xi_\la(I_\la(\psi_\la(\tilde{\ga}(t,s),\theta)))\equiv 1$ for $\theta \in [0,T(t, s))$. If $\psi_\la(\tilde{\ga}(t,s),T(t, s))\notin X^\delta$, then there exist $0<\theta_{t,s}^1<\theta_{t,s}^2<T(t, s)$ as above. Then we can prove that $I_\la(\psi_\la(\tilde{\ga}(t,s),\theta_{t,s}^2))\leq c_\la-\al_0$ as above, which contradicts the definition of $T(t, s)$. Therefore, $\ga(t,s):=\psi_\la(\tilde{\ga}(t,s), T(t,s))\in X^\delta$. Then there exists $(u,v)\in X$ such that $\|\ga(t,s)-(u,v)\|\leq \delta\leq \mathcal{C}/2$. By Lemma \ref{lemma2-1}, we have
$$\|\ga(t,s)\|\leq \|(u,v)\|+\mathcal{C}/2\leq 2C_2+\mathcal{C}.$$

To prove the continuity of $T(t, s)$, we fix any $(\tilde{t},\tilde{s})\in Q$. Assume that
$I_\la(\ga(\tilde{t},\tilde{s}))< c_\la-\al_0$ first. Then $T(\tilde{t},\tilde{s})=0$ from the definition of $T(t,s)$. So $I_\la(\tilde{\ga}(\tilde{t},\tilde{s}))< c_\la-\al_0$. By the continuity of $\tilde{\ga}$, there exists $\tau>0$ such that for any $(t,s)\in (\tilde{t}-\tau,\tilde{t}+\tau)\times (\tilde{s}-\tau,\tilde{s}+\tau)\cap Q$, we have $I_\la(\tilde{\ga}(t,s))< c_\la-\al_0$, that is, $T(t,s)=0$, and $T$ is continuous at $(\tilde{t},\tilde{s})$.

Now we assume that $I_\la(\ga(\tilde{t},\tilde{s}))= c_\la-\al_0$.
Then we have $I_\la(\tilde{\ga}(t,s))\geq c_\la-\al_0$, and thus $\tilde{\ga}(t,s)\in X^{\delta/2}$. If $\ga(\tilde{t},\tilde{s})\notin X^\delta$, then we have $T(\tilde{t},\tilde{s})>0$. From the previous proof, we can get a contradiction with the definition of $T(\tilde{t},\tilde{s})$. Therefore, we have that $\ga(\tilde{t},\tilde{s})=\psi_\la(\tilde{\ga}(\tilde{t},\tilde{s}), T(\tilde{t},\tilde{s}))\in X^\delta$, and so
$$\|I_\la'(\psi_\la(\tilde{\ga}(\tilde{t},\tilde{s}), T(\tilde{t},\tilde{s})))\|\geq l(\la)>0.$$
Then for any $\omega>0$, we have
$$I_\la(\psi_\la(\tilde{\ga}(\tilde{t},\tilde{s}), T(\tilde{t},\tilde{s})+\omega))<c_\la-\al_0.$$
By the continuity of $\psi_\la$, there exists $\tau=\tau(\omega)>0$ such that for any $(t,s)\in (\tilde{t}-\tau,\tilde{t}+\tau)\times (\tilde{s}-\tau,\tilde{s}+\tau)\cap Q$, we have $I_\la(\psi_\la(\tilde{\ga}(t,s), T(\tilde{t},\tilde{s})+\omega)+\omega))<c_\la-\al_0$, so $T(t,s)\leq T(\tilde{t},\tilde{s})+\omega$. It follows that
$$0\leq \limsup\limits_{(t,s)\rightarrow (\tilde{t},\tilde{s})}T(t,s)\leq T(\tilde{t},\tilde{s}).$$
If $T(\tilde{t},\tilde{s})=0$, we have
$$\lim\limits_{(t,s)\rightarrow (\tilde{t},\tilde{s})}T(t,s)= T(\tilde{t},\tilde{s}).$$
If $T(\tilde{t},\tilde{s})>0$, then for any $0<\omega<T(\tilde{t},\tilde{s})$, similarly we have $I_\la(\psi_\la(\tilde{\ga}(\tilde{t},\tilde{s}), T(\tilde{t},\tilde{s})+\omega)-\omega))>c_\la-\al_0$. By the continuity of $\psi_\la$ again, we have
$$\liminf\limits_{(t,s)\rightarrow (\tilde{t},\tilde{s})}T(t,s)\geq T(\tilde{t},\tilde{s}).$$
So $T$ is continuous at $(\tilde{t},\tilde{s})$. This completes the proof of {Step 2.}

Now, we have proved that $\ga(t,s)\in \Gamma$ and $\max\limits_
{(t,s)\in Q}I_\la(\ga(t,s))\leq c_\la-\al_0$, which
contradicts the definition of $c_\la$. This completes the proof. \qed

\vskip0.2in
Let's recall a version of Br\'{e}zis-Kato lemma, as in \cite[Lemma 2.5]{ZCd}.

\bl\label{lemma2-6-0} Let $\varOmega\subset \RN$ and $h\in L^\frac{N}{p}(\RN) (1<p<N)$ be a nonnegative function. Then for every $\mu>0$, there exists a constant $\sigma(\mu,h)>0$ such that
$$\int_{\varOmega}h(x)|u|^p\leq \mu\int_{\varOmega}|\nabla u|^p+\sigma(\mu,h)\int_{\varOmega}|u|^p$$
for all $u\in W^{1,p}(\varOmega)$.
\el

\vskip0.2in
\bl\label{lemma2-6} Assume that $(u,v)$ is a nontrivial solution of \eqref{1-1}, then $u,v\in L^\iy(\RN)\cap C_{loc}^{1,\beta}$ for some $\beta\in (0,1)$.
\el

\noindent {\bf Proof. } It's easily seen from \eqref{1-1} that $u\not\equiv 0,\;v\not\equiv 0$. Firstly we show that $u,v\in L^\iy(\RN)$. Set $u^+=\max\{u,0\},\;u^-=\max\{-u,0\}$, and let $A_k=\{x\in \RN:u^+(x)\leq k\},\;B_k=\RN\setminus A_k,\;k>0$.
Define
$$v_k=\begin{cases}(u^+)^{ps+1},\;\;&\mbox{in}~ A_k,\\
k^{ps}u^+,\;&\mbox{in}~ B_k,\end{cases}$$
$$w_k=\begin{cases}(u^+)^{s+1},\;\;&\mbox{in}~ A_k,\\
k^{s}u^+,\;&\mbox{in}~ B_k,\end{cases}$$
where $s\geq 0$. Then $v_k,w_k\in W^{1,p}(\RN)$, and using $v_k$ as a test function in \eqref{1-1}, we have
\be\label{2-17}
(ps+1)\int_{A_k}(u^+)^{ps}|\nabla u^+|^p+k^{ps}\int_{B_k}|\nabla u^+|^p\leq \int_{\RN}f(u)v_k+\la\int_{\RN}vv_k.
\ee
By $(F1)-(F2)$, there exists $C>0$ such that $f(t)\leq Ct^{p-1}+t^{p^*-1}$ for all $t\geq 0$. Hence
\be\label{2-18}
(ps+1)\int_{A_k}(u^+)^{ps}|\nabla u^+|^p+k^{ps}\int_{B_k}|\nabla u^+|^p\leq \int_{\RN}(C+(u^+)^{p^*-p})(u^+)^{p-1}v_k+\la\int_{\RN}vv_k.
\ee
Since
\be\label{2-19}
\int_{\RN}|\nabla w_k|^p=(s+1)^p\int_{A_k}(u^+)^{ps}|\nabla u^+|^p+k^{ps}\int_{B_k}|\nabla u^+|^p,
\ee
we see that
\be\label{2-20}
\frac{ps+1}{(s+1)^p}\int_{\RN}|\nabla w_k|^p\leq \int_{\RN}(C+(u^+)^{p^*-p})w_k^p+\la\int_{\RN}vv_k.
\ee
By Lemma \ref{lemma2-6-0}, for any $\mu > 0$, there exists $\sigma(\mu,u)$ such that
$$\int_{\RN}(u^+)^{p^*-p}w_k^p\leq \mu \int_{\RN}|\nabla w_k|^p+\sigma(\mu,u)\int_{\RN}|w_k|^p.$$
Choosing $\mu=\frac{ps+1}{2(s+1)^p}$, we have
\be\label{2-21}
\int_{\RN}|\nabla w_k|^p\leq C_s\int_{\RN}w_k^p+\la\int_{\RN}vv_k,
\ee
where $C_s=\frac{2(s+1)^p}{ps+1}(C+\sigma(\mu,u))$. By the Sobolev embedding theorem, we have
$$\big(\int_{A_k}w_k^{p^*}\big)^{p/p^*}\leq S^{-1}C_s\int_{\RN}w_k^p+S^{-1}\la\int_{\RN}vv_k,$$
where $S$ is the Sobolev best constant.

When $\frac{2N}{N+2}<p\leq \frac{2N}{N+1}$, we have
\be\label{2-22}
\begin{split}
\int_{\RN}|v|v_k&\leq \int_{\RN}|v|(u^+)^{ps+1}\\
&\leq \big(\int_{\RN}|v|^\frac{p^*s+2}{(p^*-p)s+1}\big)^\frac{(p^*-p)s+1}{p^*s+2}\big(\int_{\RN}(u^+)^{p^*s+2}\big)^{\frac{ps+1}{p^*s+2}}.
\end{split}
\ee
Therefore, we have
\be\label{2-23}
\begin{split}
\big(\int_{A_k}(u^+)^{p^*(s+1)}\big)^\frac{p}{p^*}=  &\big(\int_{A_k}|w_k|^{p^*}\big)^\frac{p}{p^*}\\
\leq &S^{-1}C_s\int_{\RN}w_k^p+S^{-1}\la \big(\int_{\RN}|v|^\frac{p^*s+2}{(p^*-p)s+1}\big)^\frac{(p^*-p)s+1}{p^*s+2}\big(\int_{\RN}(u^+)^{p^*s+2}\big)^{\frac{ps+1}{p^*s+2}}.
\end{split}
\ee
Letting $k\rightarrow \iy$, we get
\be\label{2-24}
\|u^+\|_{L^{p^*(s+1)}(\RN)}\leq \bar{C}_s\|u^+\|_{L^{p(s+1)}(\RN)}+\hat{C}_s\|v\|_{L^\frac{p^*s+2}{(p^*-p)s+1}(\RN)}^\frac{1}{p(s+1)}\|u^+\|_{L^{p^*s+2}(\RN)}^\frac{{ps+1}}{p(s+1)}.
\ee
Set $\eta_i=p^*(s_i+1),\;\xi_i=p(s_i+1),\;\zeta_i=p^*s_i+2,\;\delta_i=\frac{p^*s_i+2}{(p^*-p)s_i+1},\;i\geq 0$. Then we have $\xi_i<\zeta_i<\eta_i$ and $\delta_i\in [p,\zeta_i]$. Let $s_0=0$. Then we have $\zeta_0=2\in [p,p^*)$ and
$$\|u^+\|_{L^{p^*}(\RN)}\leq \bar{C}_0\|u^+\|_{L^p(\RN)}+\hat{C}_0\|v\|_{L^2(\RN)}^\frac{1}{p}\|u^+\|_{L^2(\RN)}^\frac{1}{p}.$$
Choosing $s_i$ such that $\zeta_{i+1}=p^*s_{i+1}+2=p^*(s_i+1)=\eta_i$, we can easily check that $\xi_{i+1},\;\delta_{i+1}\in [p,\eta_i)$ and $s_i$ is strictly increasing and tends to $+\iy$. Similar estimates for $v$ can be obtained.
Therefore, by a bootstrap argument, there exists $C=C(\|u\|,\|v\|)>0$ such that
$$\|u^+\|_{L^{(p^*-1)N}(\RN)}\leq C.$$
Moreover, by H\"{o}lder inequality, we have
$$\|u^+\|_{L^{(p-1)N}(B_2(x))}\leq C,\;\forall x\in \RN.$$
Similarly, the above conclusions also hold for $u^-,\;v^+,\;v^-$. Then we have
$$\|v\|_{L^N(B_2(x))}\leq C,\;\forall x\in \RN.$$
Now, by \cite{S}, for any ball $B_r(x)$ of radius $r$ centered at any $x\in \RN$, the solution $w\in W^{1,p}(\RN)$ of the equation $-\Delta_pw=h(x)$ satisfies the estimates
$$\sup\limits_{y\in B_1(x)}|w(y)|\leq C(N)\big(\|w\|_{L^p(B_2(x))}+\|h\|_{L^N(B_2(x))}\big).$$
By $(F1)-(F2)$, it follows that
\be\label{2-25}
\begin{split}
\sup\limits_{y\in B_1(x)}|u(y)|&\leq C(N)\big(\|u\|_{L^p(B_2(x))}+\|f(u)+\la v-|u|^{p-2}u\|_{L^N(B_2(x))}\big)\\
&\leq C(N)\big(\|u\|_{L^p(B_2(x))}+C\big),
\end{split}
\ee
for any $x\in \RN$, which implies that $u\in L^\iy(\RN)$. Similarly, we can get that $v\in L^\iy(\RN)$.

When $\frac{2N}{N+1}<p\leq 2$, we have
\be\label{2-26}
\begin{split}
	\int_{\RN}|v|v_k&\leq \int_{\RN}|v|(u^+)^{ps+1}\\
	&\leq \big(\int_{\RN}|v|^\frac{p(s+1)}{p-1}\big)^\frac{p-1}{p(s+1)}\big(\int_{\RN}(u^+)^{p(s+1)}\big)^{\frac{ps+1}{p(s+1)}}.
\end{split}
\ee
Therefore, we have
\be\label{2-27}
\begin{split}
	\big(\int_{A_k}(u^+)^{p^*(s+1)}\big)^\frac{p}{p^*}=  &\big(\int_{A_k}|w_k|^{p^*}\big)^\frac{p}{p^*}\\
	\leq &S^{-1}C_s\int_{\RN}w_k^p+S^{-1}\la \big(\int_{\RN}|v|^\frac{p(s+1)}{p-1}\big)^\frac{p-1}{p(s+1)}\big(\int_{\RN}(u^+)^{p(s+1)}\big)^{\frac{ps+1}{p(s+1)}}.
\end{split}
\ee
Letting $k\rightarrow \iy$, we get
\be\label{2-28}
\|u^+\|_{L^{p^*(s+1)}(\RN)}\leq \bar{C}_s\|u^+\|_{L^{p(s+1)}(\RN)}+\hat{C}_s\|v\|_{L^\frac{p(s+1)}{p-1}(\RN)}^\frac{1}{p(s+1)}\|u^+\|_{L^{p(s+1)}(\RN)}^\frac{{ps+1}}{p(s+1)}.
\ee
Set $\bar{\eta}_i=p^*(s_i+1),\;\bar{\xi}_i=p(s_i+1),\;\bar{\zeta}_i=\frac{p(s_i+1)}{p-1},\;i\geq 0$ and $s_0=0$. Then we have $\bar{\xi}_i\leq \bar{\zeta}_i<\bar{\eta}_i$, $\bar{\zeta}_0=\frac{p}{p-1}\in [p,p^*)$ and
$$\|u^+\|_{L^{p^*}(\RN)}\leq \bar{C}_0\|u^+\|_{L^p(\RN)}+\hat{C}_0\|u^+\|_{L^p(\RN)}^\frac{1}{p}\|v\|_{L^\frac{p}{p-1}(\RN)}^\frac{1}{p}.$$
Choosing $s_i$ such that $\bar{\zeta}_{i+1}=\frac{p(s_{i+1}+1)}{p-1}=p^*(s_i+1)=\bar{\eta}_i$, we can easily check that $\bar{\xi}_{i+1}\in [p,\bar{\eta}_i]$ and $s_i$ is strictly increasing and tends to $+\iy$.
Therefore, by a bootstrap argument, there exists $\tilde{C}=\tilde{C}(\|u\|,\|v\|)>0$ such that
$$\|u^+\|_{L^{(p^*-1)N}(\RN)}\leq \tilde{C}.$$
Moreover, by H\"{o}lder inequality, we have
$$\|u^+\|_{L^{(p-1)N}(B_2(x))}\leq \tilde{C},\;\forall x\in \RN.$$
Similarly, the above conclusions also hold for $u^-,\;v^+,\;v^-$. Then we have
$$\|v\|_{L^N(B_2(x))}\leq \tilde{C},\;\forall x\in \RN.$$
As in \eqref{2-25} that $u\in L^\iy(\RN)$. Similarly, we can get that $v\in L^\iy(\RN)$.

Thus, by \cite{T}, we know that $u,v\in C_{loc}^{1,\beta}(\RN)$ for some $\beta\in (0,1)$.

\qed

\vskip0.2in
\noindent {\bf Proof of Theorem \ref{th1}. } First we fix any $\la\in (0,\la_0)$. By Lemma \ref{lemma2-5}, there exists $\{(u_n,v_n)\}\subset X^\delta\cap I_\la^{d_\la}$ such that
$$I_\la'(u_n,v_n)\rightarrow 0~~\hbox{in}~~W_r~~\hbox{as}~ n\rightarrow \iy.$$
By Lemma \ref{lemma2-1}, $\{(u_n,v_n)\}$ are uniformly bounded in $W_r$. Up to a subsequence, we can assume that $(u_n,v_n)\rightarrow (u_\la,v_\la)$ weakly in $W_r$ and strongly in $L^{q_1}(\R^N)\times L^{q_2}(\R^N),q_1,q_2\in (p,p^*)$. Then as in the proof of Lemma \ref{lemma2-3}, we have $I_\la'(u_\la,v_\la)=0$ and $(u_\la,v_\la)$ is a solution of \eqref{1-1}. Moreover, by $(F1)-(F2)$ we have
$$\lim\limits_{n\rightarrow \iy}\int_{\R^N}(f(u_n)u_n+g(v_n)v_n)dx=\int_{\R^N}(f(u_\la)u_\la+g(v_\la)v_\la)dx.$$
By $I_\la'(u_n,v_n)(u_n,v_n)=o(1)$, we have
\be\label{2-29}
\begin{split}
&\int_{\R^N}(|\nabla u_n|^p+|\nabla v_n|^p+|u_n|^p+|v_n|^p-2\la u_nv_n)dx\\
=&\int_{\R^N}(f(u_n)u_n+g(v_n)v_n)dx+o(1)\\
=&\int_{\R^N}(f(u_\la)u_\la+g(v_\la)v_\la)dx+o(1)\\
=&\int_{\R^N}(|\nabla u_\la|^p+|\nabla v_\la|^p+|u_\la|^p+|v_\la|^p-2\la u_\la v_\la)dx+o(1).
\end{split}
\ee
Moreover, $\int_{\R^N}(u_nv_n-u_\la v_\la)=\int_{\R^N}(u_n-u_\la) v_n+\int_{\R^N}u_\la(v_n- v_\la)$ and by H\"{o}lder's inequality, it's easy to see that $\int_{\R^N}u_\la(v_n- v_\la)\rightarrow 0$. Since $p\in (\frac{2N}{N+2},2)$, we have $p<\frac{p}{p-1}, p^*>\frac{p^*}{p^*-1}$. Then there exists some $q\in (p,p^*)\cap [\frac{p^*}{p^*-1},\frac{p}{p-1}]$ such that
\be\label{2-30}
|\int_{\R^N}(u_n-u_\la) v_n|\leq (\int_{\R^N}|u_n-u_\la|^q)^\frac{1}{q}(\int_{\R^N}|v_n|^\frac{q}{q-1})^\frac{q-1}{q}\rightarrow \rightarrow 0.
\ee
As in Lemma \ref{lemma2-3}, $\nabla u_n\rightarrow \nabla u_\la, \nabla v_n\rightarrow \nabla v_\la$ a.e. $x\in \RN$. Therefore, by Br\'{e}zis-Lieb lemma and the lower semicontinuity of the norm, we have $(u_n,v_n)\rightarrow (u_\la,v_\la)$ strongly in $W_r$, and so $(u_\la,v_\la)\in X^\delta$, which implies that $u_\la\not\equiv 0,~v_\la\not\equiv 0$. Moreover, $I_\la(u_\la,v_\la)\leq d_\la$.

Let $\la_n\in (0,\la_0),n\in \mathbb{N}$, be any sequence with $\la_n\rightarrow 0$. Then by repeating the proof of Lemma \ref{lemma2-3} and passing to a subsequence, $(u_{\la_n},v_{\la_n})\rightarrow (U,V)$ strongly in $W_r$, where $U\in S_1,\;V\in S_2$. That is, $U$ is a positive radial ground state of \eqref{1-5}, and $V$ is a positive radial ground state of \eqref{1-6}. This completes the proof. \qed

\vskip0.3in
\s{Proof of Theorem \ref{th2}}
\renewcommand{\theequation}{3.\arabic{equation}}

To establish the Poho\v{z}aev's type identity for \eqref{1-1}, first we recall the Pucci-Serrin variational identity for locally Lipschitz continuous solutions of a general class of equations, see \cite[Lemma 1]{DMS} \cite[Lemma 2.14]{JS}. Let $\phi \in L_{loc}^\iy(\RN)$ and $L(s,\xi):\R\times \RN\rightarrow \R$ be a function of class $C^1$ in $s$ and $\xi$ such that for any $s\in \R$, the map $\xi\mapsto L(s,\xi)$ is strictly convex.

\vskip0.2in
\bl\label{lemma3-1-0} Let $u:\RN\rightarrow \R$ be a locally Lipschitz continuous solutions of
$$-\hbox{div}(L_\xi(u,Du))+L_s(u,Du)=\phi~~\mbox{in}~\mathcal{D}'(\RN).$$
Then
\be\label{3-1-1}
\sum_{i,j=1}^{N}\int_{\R^N}D_ih^jD_{\xi_i}L(u,Du)D_ju-\int_{\R^N}(\hbox{div} ~h)L(u,Du)=\int_{\R^N}(h\cdot Du)\phi
\ee
for every $h\in C_c^1(\RN,\RN)$.
\el

\vskip0.2in
\bl\label{lemma3-1} Assume that $(F1)-(F3)$ hold. Let $(u,v)\in W$ be a weak solution to
problem \eqref{1-1}, then we have the following Poho\v{z}aev type identity:
\be\label{3-1}
\int_{\R^N}\big(|\nabla u|^p+|\nabla v|^p\big)dx=p^*\int_{\R^N}\big(F(u)+G(v)+\la uv-\frac{1}{p}|u|^p-\frac{1}{p}|v|^p\big)dx.
\ee
\el

\noindent {\bf Proof. } By Lemma \ref{lemma2-6}, $u,v\in L^\iy(\RN)\cap C_{loc}^{1,\beta}$ for some $\beta\in (0,1)$. Let $L(u,Du)=\frac{1}{p}|Du|^p+\frac{1}{p}|u|^p-F(u)$, $\phi=v$ and $h_k(x)=T(\frac{x}{k})x$ for all $x\in \RN$ and $k\geq 1$, where $T\in C_c^1(\RN)$ satisfying $T(x)=1$ if $|x|\leq 1$ and $T(x)=0$ if $|x|\geq 2$. Then for every $k\geq 1$, we have that $h_k\in C_c^1(\RN,\RN)$ and
$$D_ih_k^j(x)=D_iT(\frac{x}{k})\frac{x_j}{k}+T(\frac{x}{k})\delta_{ij},$$
$$(\hbox{div}~h_k)(x)=DT(\frac{x}{k})\cdot \frac{x}{k}+NT(\frac{x}{k}).$$
Then by \eqref{3-1-1}, we have
\be\label{3-1-2}
\begin{split}
\int_{\R^N}(T(\frac{x}{k})x\cdot Du)v=&\sum_{i,j=1}^{N}\int_{\R^N}D_iT(\frac{x}{k})\frac{x_j}{k}D_{\xi_i}L(u,Du)D_ju+\int_{\R^N}T(\frac{x}{k})D_\xi L(u,Du)\cdot Du\\
&-\int_{\R^N}DT(\frac{x}{k})\cdot \frac{x}{k}L(u,Du)-\int_{\R^N}NT(\frac{x}{k})L(u,Du).
\end{split}
\ee
Similarly, we have
\be\label{3-1-3}
\begin{split}
	\int_{\R^N}(T(\frac{x}{k})x\cdot Dv)u=&\sum_{i,j=1}^{N}\int_{\R^N}D_iT(\frac{x}{k})\frac{x_j}{k}D_{\xi_i}\bar{L}(v,Dv)D_jv+\int_{\R^N}T(\frac{x}{k})D_\xi \bar{L}(v,Dv)\cdot Dv\\
	&-\int_{\R^N}DT(\frac{x}{k})\cdot \frac{x}{k}\bar{L}(v,Dv)-\int_{\R^N}NT(\frac{x}{k})\bar{L}(v,Dv),
\end{split}
\ee
where $\bar{L}(v,Dv)=\frac{1}{p}|Dv|^p+\frac{1}{p}|v|^p-G(v)$. Integrating by parts, we have
\be\label{3-1-4}
\begin{split}
\int_{\R^N}(T(\frac{x}{k})x\cdot Du)v=&\int_{B_{2k}(0)}(T(\frac{x}{k})x\cdot Du)v\\
=&-\int_{B_{2k}(0)}(T(\frac{x}{k})x\cdot Dv)u-\int_{B_{2k}(0)}uv(NT(\frac{x}{k})+DT(\frac{x}{k})\cdot \frac{x}{k}).
\end{split}
\ee
Combining \eqref{3-1-2} with \eqref{3-1-3}, we have
\be\label{3-1-5}
\begin{split}
-\int_{\R^N}uv(NT(\frac{x}{k})+DT(\frac{x}{k})\cdot \frac{x}{k})
=&\sum_{i,j=1}^{N}\int_{\R^N}D_iT(\frac{x}{k})\frac{x_j}{k}(D_{\xi_i}L(u,Du)D_ju+D_{\xi_i}\bar{L}(v,Dv)D_jv)\\
&+\int_{\R^N}T(\frac{x}{k})(D_\xi L(u,Du)\cdot Du+D_\xi \bar{L}(v,Dv)\cdot Dv)\\
&-\int_{\R^N}DT(\frac{x}{k})\cdot \frac{x}{k}(L(u,Du)+\bar{L}(v,Dv))\\
&-\int_{\R^N}NT(\frac{x}{k})(L(u,Du)+\bar{L}(v,Dv)).
\end{split}
\ee
Since there exists $C>0$ such that $|D_iT(\frac{x}{k})\frac{x_j}{k}|\leq C$ for every $x\in \RN,\;k\geq 1,\;\;i,j=1,...,N$, by the Dominated Convergence Theorem, we can obtain \eqref{3-1}.  \qed

\vskip0.2in
\bl\label{lemma3-2} Assume that $(F1)-(F3)$ hold, then

\begin{itemize}
	\item[(1)]  $\mathcal{M}_\la$ is a $C^1$ manifold;
	\item[(2)]  there exists a positive constant $\rho_\la > 0$ such that $\|(u,v)\|\geq \rho_\la$ for all $(u,v)\in \mathcal{M}_\la$;
	\item[(3)]  for any $(u,v)\in W\setminus\{(0,0)\}$ with $P(u, v)\leq0$, there exists a unique $t_{u,v} > 0$ such that $(u^{t_{u,v}},v^{t_{u,v}})\in \mathcal{M}_\la$, where $u^t(x)=u(x/t),\;v^t(x)=v(x/t)$. Moreover, the value $t_{u,v}$ is the maximum point of the function $t\mapsto I_\la(u^t,v^t)$. In particular, if $P(u, v) < 0$, then
	$t_{u,v}\in(0, 1)$; if $P(u, v)= 0$, then $t_{u,v} = 1$.
\end{itemize}
\el

\noindent {\bf Proof. } (1) Since $P(u,v)$ is a $C^1$ functional,
in order to prove $\mathcal{M}_\la$ is a $C^1$ manifold,
it suffices to prove that $P'(u,v)\neq 0$ for all $(u,v)\in \mathcal{M}_\la$.
Indeed, assume by contradiction that $P'(u,v)=0$ for some $(u,v)\in \mathcal{M}_\la$.
Then in a weak sense, $(u,v)$ can be seen as a solution of the problem
\be\label{3-2}
\begin{cases}-\Delta_p u +\frac{N}{N-p}|u|^{p-2}u =
	\frac{N}{N-p}f(u)+\frac{N}{N-p}\lambda v, \quad x\in \R^N,\\
	-\Delta_p v +\frac{N}{N-p}|v|^{p-2}v =
	\frac{N}{N-p}g(v)+\frac{N}{N-p}\lambda u, \quad x\in \R^N.\end{cases}
\ee
As a consequence, we see that $(u,v)$ satisfies the Poho\v{z}aev type identity
\be\label{3-3}
\int_{\R^N}\big(|\nabla u|^p+|\nabla v|^p\big)dx=\frac{{p^*}^2}{p}\int_{\R^N}\big(F(u)+G(v)+\la uv-\frac{1}{p}|u|^p-\frac{1}{p}|v|^p\big)dx.
\ee
Since $P(u, v)= 0$, we deduce that
$$(1-\frac{N}{N-p})\int_{\R^N}\big(|\nabla u|^p+|\nabla v|^p\big)dx=0,$$
which implies that $u = 0$ and $v = 0$, a contradiction since $(u,v)\in \mathcal{M}_\la$. Then $\mathcal{M}_\la$ is a $C^1$ manifold.

(2) First, by $(F1)-(F2)$, for any $\varepsilon>0$, there exists $C_\varepsilon>0$ such that
$$|F(t)|,|G(t)|\leq \varepsilon |t|^p+C_\varepsilon|t|^{p^*},\;\forall t\in \R.$$
If $(u,v)\in \mathcal{M}_\la$, we have
\be\label{3-4}
\begin{split}
	\int_{\R^N}\big(|\nabla u|^p+|\nabla v|^p\big)dx=&p^*\int_{\R^N}\big(F(u)+G(v)+\la uv\big)dx-\frac{p^*}{p}\int_{\R^N}\big(|u|^p+|v|^p\big)dx\\
	\leq &p^*C_\varepsilon \int_{\R^N}\big(|u|^{p^*}+|v|^{p^*}\big)dx+\la p^*\int_{\R^N} uv dx\\
	&-\frac{p^*}{p}(1-p\varepsilon)\int_{\R^N}\big(|u|^p+|v|^p\big)dx.
\end{split}
\ee
When $p\in [\frac{2N}{N+1},2)$, by H\"{o}lder's inequality and Young's inequality, we have
\be\label{3-4-1}
\begin{split}
\int_{\R^N} uv&\leq \int_{\R^N} \frac{\varepsilon|u|^p}{p}+\frac{|v|^q}{\varepsilon^\frac{q}{p}q}\\
&\leq \frac{\varepsilon}{p}\int_{\R^N}|u|^p +\frac{1}{\varepsilon^\frac{q}{p}q}(\int_{\R^N}|v|^p)^\theta(\int_{\R^N}|v|^{p^*})^{1-\theta}\\
&\leq \frac{\varepsilon}{p}\int_{\R^N}|u|^p +\frac{1}{q}\big(\theta\varepsilon\int_{\R^N}|v|^p+(1-\theta)\varepsilon^\frac{q/p+\theta}{\theta-1}\int_{\R^N}|v|^{p^*}\big),
\end{split}
\ee
where $q=\frac{p}{p-1}\in [p,p^*]$ and $\theta=\frac{(N+1)p-2N}{p}\in[0,1)$.
When $p\in [\frac{2N}{N+2},\frac{2N}{N+1})$, we have
\be\label{3-4-2}
\begin{split}
	\int_{\R^N} uv&\leq \int_{\R^N} \frac{|u|^{p^*}}{p^*}+\frac{|v|^q}{q}\\
	&\leq \frac{1}{p^*}\int_{\R^N}|u|^{p^*} +\frac{1}{q}(\int_{\R^N}|v|^p)^\theta(\int_{\R^N}|v|^{p^*})^{1-\theta}\\
	&\leq \frac{1}{p^*}\int_{\R^N}|u|^{p^*} +\frac{1}{q}\big(\theta\varepsilon\int_{\R^N}|v|^p+(1-\theta)\varepsilon^\frac{\theta}{\theta-1}\int_{\R^N}|v|^{p^*}\big),
\end{split}
\ee
where $q=\frac{p^*}{p^*-1}\in [p,p^*]$ and $\theta=\frac{(N+2)p-2N}{p}\in [0,1)$. Then we can choose $\varepsilon>0$ small such that
\be\label{3-4-3}
\begin{split}
\int_{\R^N}\big(|\nabla u|^p+|\nabla v|^p\big)dx
\leq &C\int_{\R^N}\big(|u|^{p^*}+|v|^{p^*}\big)dx\\
\leq &C\big(\int_{\R^N}|\nabla u|^pdx\big)^{\frac{p^*}{p}}+C\big(\int_{\R^N}|\nabla v|^pdx\big)^{\frac{p^*}{p}},
\end{split}
\ee
which implies that there exists $\rho_\la > 0$ such that
$$\|(u,v)\|\geq \big(\int_{\R^N}\big(|\nabla u|^p+|\nabla v|^p\big)dx\big)^{1/p}\geq \rho_\la.$$

(3) Let $(u,v)\in W\setminus \{(0,0)\}$ with $P(u, v)\leq0$ and define
\be\nonumber
\begin{split}
	h(t):=& I_\la(u^t,v^t)\\
	=&\frac{t^{N-p}}{p}\int_{\R^N}\big(|\nabla u|^p+|\nabla v|^p\big)dx\\
	&-t^N\int_{\R^N}\big(F(u)+G(v)+\la uv\big)dx-\frac{1}{p}|u|^p-\frac{1}{p}|v|^p\big)dx.
\end{split}
\ee
Then we obtain that $h(t)>0$ for $t>0$ small enough.Since $P(u, v)\leq0$, it's easy to see that $h(t)\rightarrow -\iy$ as $t\rightarrow +\iy$. Hence there exists $t_{u,v} > 0$ such that $h(t_{u,v})=\max\limits_{t\geq 0}h(t)$ and $h'(t_{u,v})=0$. Note that $P(u^t,v^t)=\frac{p}{N-p}th'(t)$, so we have $P(u^{t_{u,v}},v^{t_{u,v}})=0$. Moreover, if $P(u^{t_{u,v}},v^{t_{u,v}})=0$, we have
$$t_{u,v}^{N-p}\int_{\R^N}\big(|\nabla u|^p+|\nabla v|^p\big)dx=p^*t_{u,v}^N\int_{\R^N}\big(F(u)+G(v)+\la uv-\frac{1}{p}|u|^p-\frac{1}{p}|v|^p\big)dx,$$
and
\be\label{3-4-4}
t_{u,v}^p=\frac{\int_{\R^N}\big(|\nabla u|^p+|\nabla v|^p\big)dx}{p^*\int_{\R^N}\big(F(u)+G(v)+\la uv-\frac{1}{p}|u|^p-\frac{1}{p}|v|^p\big)dx}.
\ee
Thus $t_{u,v}$ is the unique critical point of $h(t)$ and the conclusions hold.\qed

\vskip0.2in
\bl\label{lemma3-3} Let $\{(u_n,v_n)\}\subset \mathcal{M}_\la$ be a bounded sequence. Then there exist a sequence $\{y_n\}\subset \R^N$ and constants $R,\theta>0$ such that
$$\liminf\limits_{n\rightarrow \iy}\int_{B_R(y_n)}\big(|u_n|^p+|v_n|^p\big)dx\geq \theta>0.$$
\el

\noindent {\bf Proof. } Assume by contradiction that for any $R>0$, up to a subsequence, there hold
\be\label{3-5}
\lim\limits_{n\rightarrow \iy}\sup\limits_{y\in \R^N}\int_{B_R(y)}|u_n|^pdx=0,\;\;\lim\limits_{n\rightarrow \iy}\sup\limits_{y\in \R^N}\int_{B_R(y)}|v_n|^pdx=0.
\ee
Then we have
$$u_n,v_n\rw 0~~\hbox{in}~~L^q(\R^N)~~\hbox{for~all}~~q\in (p,p^*).$$
By $(F1)-(F2)$ and \eqref{2-30}, we have that $\int_{\R^N}(F(u_n)), \int_{\R^N}(G(v_n)), \int_{\R^N}u_nv_n \rightarrow 0$.
Since $\{(u_n,v_n)\}\subset \mathcal{M}_\la$, it's easy to see that $(u_n,v_n)\rw (0,0)$ in $W$, contradicting Lemma \ref{lemma3-2} (2). This completes the proof. \qed

\vskip0.2in
Now we define
\be\label{3-6}
m_\la =\inf\limits_{(u,v)\in \mathcal{M}_\la}I_\la(u,v).
\ee

\bl\label{lemma3-4} For any $\la>0$, we have $m_\la>0$.
\el

\noindent {\bf Proof. } If $(u,v)\in \mathcal{M}_\la$, we have
$$I_\la(u,v)=\frac{1}{N}\int_{\R^N}\big(|\nabla u|^p+|\nabla v|^p\big)dx\geq \frac{\rho_0}{N}.$$
It follows that $m_\la>0$.  \qed

\vskip0.2in
\bl\label{lemma3-5} If $m_\la$ is attained at $(u,v)\in \mathcal{M}_\la$, then $(u,v)$ is a solution of \eqref{1-1}.\el

\noindent {\bf Proof. }  Assume $(u,v)\in \mathcal{M}_\la$ such that $I_\la(u,v)=m_\la$. Then by the Lagrange multiplier theorem, there exists a Lagrange multiplier $\delta\in \R$ such that
$$I_\la'(u,v)=\delta P'(u,v).$$
Then in a weak sense, $(u,v)$ can be seen as a solution of the problem
\be\nonumber
\begin{cases}-(1-\delta p)\Delta_p u +(1-\delta p^\ast)|u|^{p-2}u =
	(1-\delta p^\ast)f(u)+(1-\delta p^\ast)\al\lambda |u|^{\al-2}u|v|^\beta, \quad x\in \R^N,\\
	-(1-\delta p)\Delta_p v +(1-\delta p^\ast)|v|^{p-2}v =
	(1-\delta p^\ast)g(v)+(1-\delta p^\ast)\beta\lambda |u|^\alpha|v|^{\beta-2}v, \quad x\in \R^N.\end{cases}
\ee
As a consequence, $(u,v)$ satisfies the following Poho\v{z}aev type identity
\be\label{3-7}
(1-\delta p)\int_{\R^N}\big(|\nabla u|^p+|\nabla v|^p\big)=p^\ast(1-\delta p^\ast)\int_{\R^N}\big(F(u)+G(v)+\la |u|^\al|v|^\beta-\frac{1}{p}|u|^p-\frac{1}{p}|v|^p\big).
\ee
Since $P(u,v)=0$, we get that $\delta=0$. Thus, we have $I_\la'(u,v)=0$, namely, $(u,v)$ is a solution of \eqref{1-1}.  \qed

\vskip0.2in
\bl\label{lemma3-6} For any $\la>0$ and $(u,v)\in \mathcal{M}_\la$, there exists $(\bar{u},\bar{v})\in \mathcal{M}_\la\cap W_r$ such that $I_\la(\bar{u},\bar{v})\leq I_\la(u,v)$.
\el

\noindent {\bf Proof. } Let $(u,v)\in \mathcal{M}_\la$, and $(u^*,v^*)$ be the Schwarz symmetric radial decreasing rearrangement of $(u,v)$. Then we have
$$\int_{\R^N}|\nabla u^*|^pdx\leq \int_{\R^N}|\nabla u|^pdx,\;\;\int_{\R^N}|\nabla v^*|^pdx\leq \int_{\R^N}|\nabla v|^pdx,$$
$$\int_{\R^N}|u^*|^pdx= \int_{\R^N}|u|^pdx,\;\;\int_{\R^N}|v^*|^pdx= \int_{\R^N}|v|^pdx,$$
$$\int_{\R^N}F(u^*)dx= \int_{\R^N}F(u)dx,\;\;\int_{\R^N}G(v^*)dx= \int_{\R^N}G(v)dx,\;\;\int_{\R^N}u^*v^*dx\geq \int_{\R^N}uvdx.$$
Then we obtain
\be\label{3-8}
\begin{split}
	&\int_{\R^N}\big(|\nabla u^*|^p+|\nabla v^*|^p\big)dx\\
	\leq &\int_{\R^N}\big(|\nabla u|^p+|\nabla v|^p\big)dx\\
	=&p^*\int_{\R^N}\big(F(u)+G(v)+\la uv-\frac{1}{p}|u|^p-\frac{1}{p}|v|^p\big)dx\\
	\leq &p^*\int_{\R^N}\big(F(u^*)+G(v^*)+\la u^*v^*-\frac{1}{p}|u^*|^p-\frac{1}{p}|v^*|^p\big)dx.
\end{split}
\ee
Then by Lemma \ref{lemma3-2} (3), there exists $\bar{t}\in (0,1]$ such that $(\bar{u},\bar{v}):=(u^*(\frac{\cdot}{\bar{t}}),v^*(\frac{\cdot}{\bar{t}}))\in \mathcal{M}_\la\cap W_r$. Then we have
\be\label{3-9}
\begin{split}
	I_\la(\bar{u},\bar{v})=&(\frac{1}{p}-\frac{1}{p^*})\int_{\R^N}\big(|\nabla \bar{u}|^p+|\nabla \bar{v}|^p\big)dx\\
	= &\frac{\bar{t}^{N-p}}{N}\int_{\R^N}\big(|\nabla u^*|^p+|\nabla v^*|^p\big)dx\\
	\leq &\frac{1}{N}\int_{\R^N}\big(|\nabla u^*|^p+|\nabla v^*|^p\big)dx\\
	\leq &I_\la(u,v).
\end{split}
\ee
This completes the proof.  \qed

\vskip0.2in
\bl\label{lemma3-7} For any $\la>0$, there exists $(u_\la,v_\la)\in \mathcal{M}_\la$ such that $I_\la(u_\la,v_\la)=m_\la$ and $u_\la,v_\la>0$.
\el

\noindent {\bf Proof. } Let $\{(u_n,v_n)\}$ be a minimizing sequence for $m_\la$. By Lemma \ref{lemma3-6}, we can assume that $\{(u_n,v_n)\}\subset \mathcal{M}_\la\cap W_r$ and $u_n,v_n\geq 0$. We claim that $\{(u_n,v_n)\}$ is bounded. Indeed, since $P(u_n,v_n)=0$, we have
$$I_\la(u_n,v_n)=\frac{1}{N}\int_{\R^N}\big(|\nabla u_n|^p+|\nabla v_n|^p\big)dx.$$
Then $\{u_n\},\{v_n\}$ are bounded in $\mathcal{D}^{1,p}(\R^N)$. Moreover, by \eqref{3-4} and the Sobolev embedding theorem, we deduce the boundedness of the $L^p$ norm of $\{u_n\},\{v_n\}$. Therefore $\{(u_n,v_n)\}$ is bounded in $W_r$. Up to a subsequence, we can assume $(u_n,v_n)\rw (u,v)$ weakly in $W_r$ and strongly in $L^q(\R^N),\; q\in (p,p^*)$. By Lemma \ref{lemma3-3}, we know that there exists a sequence $\{y_n\}\subset \R^N$ and constants $R,\theta>0$ such that
$$\liminf\limits_{n\rightarrow \iy}\int_{B_R(y_n)}\big(|u_n|^p+|v_n|^p\big)dx\geq \theta>0.$$
Now we show that $\{y_n\}$ is bounded. Indeed, if $\{y_n\}$ is unbounded, then there exists a subsequence $\{y_{n_k}\}$ such that $\{|y_{n_k}|\}$ is increasing and tends to $+\iy$. Without loss of generality, we assume that 
$$\liminf\limits_{n\rightarrow \iy}\int_{B_R(y_n)}|u_n|^pdx\geq \frac{1}{2}\theta>0.$$
Since  $u_n\in W_r^{1,p}(\RN)$, by Radial Lemma \cite[Lemma A.IV]{BeL}, it is easy to see
that $|u_n(x)|\leq C|x|^{-\frac{N}{p}}$, with $C$ independent of $n$. Then we have
$$\int_{B_R(y_{n_k})}|u_{n_k}|^pdx\leq C\int_{B_{|y_{n_k}|+R}(0)\setminus B_{|y_{n_k}|-R}(0)}|x|^{-N}dx=C\ln\frac{|y_{n_k}|+R}{|y_{n_k}|-R}\rightarrow 0$$
as $k\rightarrow +\iy$, which is a contradiction. Then $\{y_n\}$ is bounded and there exists $R_1>R$ such that 
$$\liminf\limits_{n\rightarrow \iy}\int_{B_{R_1}(0)}\big(|u_n|^p+|v_n|^p\big)dx\geq \theta>0.$$
Up to a subsequence, we can assume $(u_n,v_n)\rw (\bar{u},\bar{v})$ weakly in $W$ and strongly in $L^{p_1}(\R^N)\times L^{p_2}(\R^N),\; p_1,\;p_2\in (p,p^*)$.
Passing to the limit, we get that
$$\int_{B_{R_1}(0)}\big(|\bar{u}|^p+|\bar{v}|^p\big)dx\geq \theta>0,$$
which implies that $(\bar{u},\bar{v})\neq (0,0)$. Moreover, by $(F1)-(F2)$ and \cite[Lemma A.I.]{BeL}, we have that $\int_{\R^N}(F(u_n))\rightarrow \int_{\R^N}(F(\bar{u})), \int_{\R^N}(G(v_n))\rightarrow \int_{\R^N}(G(\bar{v}))$. By Fatou's lemma, we can also deduce that
\be\label{3-10-0}
\int_{\R^N}\big(\frac{\mu}{p}|u|^p+\frac{\nu}{p}|v|^p-\la |u|^\al|v|^\beta\big)dx\leq \varliminf\limits_{n\rightarrow +\iy} \int_{\R^N}\big(\frac{\mu}{p}|u_n|^p+\frac{\nu}{p}|v_n|^p-\la |u_n|^\al|v_n|^\beta\big)dx.
\ee
Then we have $P(\bar{u},\bar{v})\leq \varliminf\limits_{n\rightarrow \iy}P(\bar{u}_n,\bar{v}_n)=0$.
Therefore, by Lemma \ref{lemma3-2} (3), there exists $t\in (0,1]$ such that $(u_\la,v_\la):=(\bar{u}^t,\bar{v}^t)\in \mathcal{M}_\la$. Then we have
\be\label{3-10}
\begin{split}
	m_\la\leq J_\la(u_\la,v_\la)=&\frac{1}{N}\int_{\R^N}\big(|\nabla u_\la|^p+|\nabla v_\la|^p\big)dx\\
	\leq &\frac{1}{N}\int_{\R^N}\big(|\nabla \bar{u}|^p+|\nabla \bar{v}|^p\big)dx\\
	\leq &\varliminf\limits_{n\rightarrow \iy}\frac{1}{N}\int_{\R^N}\big(|\nabla \bar{u}_n|^p+|\nabla \bar{v}_n|^p\big)dx\\
	=&\lim\limits_{n\rightarrow \iy}I_\la(\bar{u}_n,\bar{v}_n)=m_\la.
\end{split}
\ee
Hence $(u_\la,v_\la)$ is a minimizer of $I_\la$. By Lemma \ref{lemma3-5}, $(u_\la,v_\la)$ is a solution of \eqref{1-1}. Then by the strong maximum principle, we can get that $u_\la,v_\la>0$.  \qed

In order to study the asymptotic behavior for the vector ground state solutions with
respect to the parameter $\la$, we need the following result.

\vskip0.2in
\bl\label{lemma3-8} The map $\la \mapsto m_\la,\;\la\geq 0$ is strictly decreasing.
\el

\noindent {\bf Proof. } For given $\la_1<\la_2$, let $(u_i,v_i)\in \mathcal{M}_{\la_i}$ be such that $m_{\la_i}=I_{\la_i}(u_i,v_i),u_i,v_i>0,\;i=1,2$. Choose $t>0$ such that $(u_1^t,v_1^t)\in \mathcal{M}_{\la_2}$, that is,
\be\label{3-11}
\begin{split}
	&t^{N-p}\int_{\R^N}\big(|\nabla u_1|^p+|\nabla v_1|^p\big)dx\\
	=&p^*t^N\int_{\R^N}\big(F(u_1)+G(v_1)+\la_2 u_1v_1\big)dx-\frac{1}{p}|u_1|^p-\frac{1}{p}|v_1|^p\big)dx.
\end{split}
\ee
Since $(u_1,v_1)\in \mathcal{M}_{\la_1}$, we have
\be\label{3-12}
\begin{split}
	&\int_{\R^N}\big(|\nabla u_1|^p+|\nabla v_1|^p\big)dx\\
	=&p^*\int_{\R^N}\big(F(u_1)+G(v_1)+\la_1 u_1v_1\big)dx-\frac{1}{p}|u_1|^p-\frac{1}{p}|v_1|^p\big)dx\\
	<&p^*\int_{\R^N}\big(F(u_1)+G(v_1)+\la_2 u_1v_1\big)dx-\frac{1}{p}|u_1|^p-\frac{1}{p}|v_1|^p\big)dx.
\end{split}
\ee
Then we deduce that $t<1$ and we have
\be\label{3-13}
\begin{split}
	m_{\la_2}\leq &I_{\la_2}(u_1^t,v_1^t)\\
	=&\frac{t^{N-p}}{N}\int_{\R^N}\big(|\nabla u_1|^p+|\nabla v_1|^p\big)dx\\
	<&\frac{1}{N}\int_{\R^N}\big(|\nabla u_1|^p+|\nabla v_1|^p\big)dx\\
	=&I_{\la_1}(u_1,v_1)=m_{\la_1}.
\end{split}
\ee
\qed

\vskip0.2in
\noindent {\bf Proof of Theorem \ref{th2}. } By Lemma \ref{lemma3-7}, for any $\la>0$, \eqref{1-1} has a positive radial ground state $(u_\la,v_\la)$. Let $\{\la_n\}\subset (0,+\iy)$ be a sequence with $\la_n\searrow 0$ as $n\rightarrow \iy$ (we assume that $\la_n<1$) and $\{(u_{\la_n},v_{\la_n})\}\subset W_r$ be a sequence of positive vector ground state solutions. Then we have
\be\label{3-17}
I_{\la_n}(u_{\la_n},v_{\la_n})=m_{\la_n},\;I_{\la_n}'(u_{\la_n},v_{\la_n})=0,\;P_{\la_n}(u_{\la_n},v_{\la_n})=0.
\ee
As in Lemma\ref{lemma3-7}, we know that $\{(u_{\la_n},v_{\la_n})\}$ is bounded in $W_r$. Up to a subsequence we may assume that
\be\label{3-18}
\begin{cases}(u_{\la_n},v_{\la_n})\rightharpoonup (u_0,v_0), \quad \mbox{in}~ W_r,\\
	(u_{\la_n},v_{\la_n})\rightarrow (u_0,v_0), \quad \mbox{for}~ a.e. x\in \R^N,\\
	(u_{\la_n},v_{\la_n})\rightarrow (u_0,v_0), \quad \mbox{in}~ L^{q_1}(\RN)\times L^{q_2}(\RN), p<q_1,q_2<p^*.
\end{cases}
\ee
Then $u_0,v_0\geq 0$ and are radial. By \eqref{2-30} and \eqref{3-17}, for $n>0$ large, we have
\be\label{3-19}
\begin{split}
	 \int_{\R^N}\big(F(u_{\la_n})+G(v_{\la_n}))\big)dx&=\frac{p}{p^*-p}m_{\la_n}+\int_{\R^N}\big(\frac{1}{p}|u_{\la_n}|^p+\frac{1}{p}|v_{\la_n}|^p-\la_nu_{\la_n}v_{\la_n}\big)dx\\
	&\geq \frac{p}{2(p^*-p)}m_{1}>0.
\end{split}
\ee
On the other hand, by \cite[Theorem A.I.]{BeL}, we obtain
$$\lim\limits_{n\rightarrow \iy}\int_{\R^N}F(u_{\la_n})dx=\int_{\R^N}F(u_0)dx,\;\;\; \lim\limits_{n\rightarrow \iy}\int_{\R^N}G(v_{\la_n})dx=\int_{\R^N}G(v_0)dx.$$
Then we conclude that $(u_0,v_0)\neq (0,0)$.

Since $I_0'(u_{\la_n},v_{\la_n})\rightarrow 0$ and $\lim\limits_{n\rightarrow \iy}I_0(u_{\la_n},v_{\la_n})=\lim\limits_{n\rightarrow \iy}I_{\la_n}(u_{\la_n},v_{\la_n})\in (m_1,m_0]$, as in the proof of Lemma \ref{lemma2-3}, we have that $I_0'(u_0,v_0)=0$ and
\be\label{3-21}
\begin{split}
0< I_0(u_0,v_0)=&\frac{p}{N}\int_{\R^N}\big(|\nabla u_0|^p+|\nabla v_0|^p\big)dx\\
\leq&\liminf\limits_{n\rightarrow \iy}\frac{p}{N}\int_{\R^N}\big(|\nabla u_{\la_n}|^p+|\nabla v_{\la_n}|^p\big)dx\\
=&\lim\limits_{n\rightarrow \iy}I_{\la_n}(u_{\la_n},v_{\la_n})=\lim\limits_{n\rightarrow \iy}m_{\la_n}\leq m_0.
\end{split}
\ee
Then we obtain $I_0(u_0,v_0)=m_0$ and $(u_0,v_0)$ is a ground state of \eqref{1-1} for $\la=0$. Furthermore, by \eqref{3-21},
$$\int_{\R^N}\big(|\nabla u_{\la_n}|^p+|\nabla v_{\la_n}|^p\big)dx\rightarrow \int_{\R^N}\big(|\nabla u_0|^p+|\nabla v_0|^p\big)dx$$
and
$$\int_{\R^N}\big(F(u_{\la_n})+G(v_{\la_n})+\la_n u_{\la_n}v_{\la_n}-\frac{1}{p}|u_{\la_n}|^p-\frac{1}{p}|v_{\la_n}|^p\big)dx\rightarrow \int_{\R^N}\big(F(u_0)+G(v_0)-\frac{1}{p}|u_0|^p-\frac{1}{p}|v_0|^p\big)dx,$$
then we deduce that
$$\int_{\R^N}\big(|u_{\la_n}|^p+|v_{\la_n}|^p\big)dx\rightarrow \int_{\R^N}\big(|u_0|^p+|v_0|^p\big)dx.$$
Therefore, we have $(u_{\la_n},v_{\la_n})\rightarrow (u_0,v_0)$ strongly in $W_r$. Since $(u_0,v_0)$ is a ground state of \eqref{1-1} with $\la=0$, if $u_0\neq 0, v_0\neq 0$, then we have $I_0(u_0,v_0)>I_0(0,v_0)$ and $I_0(u_0,v_0)>I_0(u_0,0)$.
Therefore we have either $(u_0,v_0)=(u_0,0)$ or $(u_0,v_0)=(0,v_0)$.
\qed

\end{document}